\documentclass[12pt]{article}

\usepackage{epsfig}
\usepackage{geometry}
\usepackage{amsfonts}
\usepackage{amsmath}
\usepackage{enumitem}
\usepackage{bbm}

\begin{document}

\setcounter{page}{1} \setcounter{section}{0}
\newtheorem{theorem}{Theorem}[section]
\newtheorem{lemma}[theorem]{Lemma}
\newtheorem{corollary}[theorem]{Corollary}
\newtheorem{proposition}[theorem]{Proposition}
\newtheorem{observation}[theorem]{Observation}
\newtheorem{definition}[theorem]{Definition}
\newtheorem{claim}{Claim}
\newtheorem{conjecture}[theorem]{Conjecture}
\newtheorem{problem}[theorem]{Problem}

\title{Packing Steiner Trees}

\author{
   Matt DeVos\thanks{Department of Mathematics, Simon Fraser University, Burnaby, B.C., Canada V5A 1S6, mdevos@sfu.ca.
     Supported in part by an NSERC Discovery Grant (Canada) and a Sloan Fellowship.}
\and
   Jessica McDonald\thanks{Department of Mathematics and Statistics, Auburn University, Auburn, AL, USA 36849, mcdonald@auburn.edu.}
\and
   Irene Pivotto\thanks{School of Mathematics and Statistics, University of Western Australia, Perth, WA, Australia 6009, irene.pivotto@uwa.edu.au.}
}

\date{}

\maketitle

\begin{abstract}
Let $T$ be a distinguished subset of vertices in a graph $G$.  A $T$-\emph{Steiner tree} is a subgraph of $G$ that
is a tree and that spans $T$. Kriesell conjectured that $G$ contains $k$ pairwise edge-disjoint $T$-Steiner trees provided that every edge-cut of $G$ that separates $T$ has size $\ge 2k$.
When $T=V(G)$ a $T$-Steiner tree is a spanning tree and the conjecture is a consequence of a classic theorem due to Nash-Williams and Tutte. Lau proved that Kriesell's conjecture holds when $2k$ is replaced by $24k$, and recently West and Wu have lowered this value to $6.5k$.
Our main result makes a further improvement to $5k+4$.
\end{abstract}

\section{Introduction}

Graphs in this paper are permitted to have multiple edges and loops (although loops will have no effect on our problem).  A classic theorem due to Nash-Williams and Tutte gives a necessary and sufficient condition for a graph to have $k$ pairwise edge-disjoint spanning trees.  Before stating their theorem, let us introduce some notation. Given a graph or hypergraph $H$ with $F\subseteq E(H)$ and ${\mathcal P}$ a partition of $V(H)$, we let
$$\Lambda^{in}_{\mathcal P}(F) = \{ e \in F \mid \mbox{$e$ is contained in a block of ${\mathcal P}$} \}$$
and $\Lambda^{out}_{\mathcal P}(F) = F \setminus \Lambda^{in}_{\mathcal P}(F),$ with $\lambda^{in}_{\mathcal P}(F) = | \Lambda^{in}_{\mathcal P}(F) |$ and $\lambda^{out}_{\mathcal P}(F) = | \Lambda^{out}_{\mathcal P}(F) |$.

\begin{theorem}[Nash-Williams \cite{NW}; Tutte \cite{Tu}]\label{nwt_theorem}
A graph $G = (V,E)$ has $k$ pairwise edge disjoint spanning trees if and only if every partition ${\mathcal P}$ of $V$ satisfies
\begin{equation}\label{lambdaout}
 \lambda^{out}_{\mathcal P}(E) \ge k ( |{\mathcal P}| - 1 ).
 \end{equation}
\end{theorem}

Let us remark that the ``only if'' direction of this theorem is obvious since each spanning tree must contain at least $| {\mathcal P} | - 1$ edges from $\Lambda^{out}_{\mathcal P}(E)$.  This result has a natural generalization to matroids and this is, in some sense, why we have such a nice characterization for it.  One immediate and appealing corollary of this result is the following.

\begin{corollary}[Nash-Williams \cite{NW}; Tutte \cite{Tu}]
\label{nwt_cor}
Every $2k$-edge-connected graph contains $k$ pairwise edge-disjoint spanning trees.
\end{corollary}

Let $G = (V,E)$ be a graph and fix a subset $T \subseteq V$ of vertices.  A $T$-\emph{Steiner tree} is a subgraph of $G$ which is a tree that spans $T$.  In this paper we are interested in packing (i.e. finding edge disjoint) $T$-Steiner trees.  If $G$ is a graph with a nontrivial edge cut of size less than $k$, then $G$ cannot  contain $k$ edge-disjoint spanning trees (since every spanning tree meets every nontrivial edge cut);  However, such a graph may still contain $k$ edge-disjoint $T$-Steiner trees, so long as this edge cut does not have vertices of $T$ on either side.  Assuming $G$ is connected, we say that an edge-cut $C$ \emph{separates} $T$ if there are at least two components of $G \setminus C$ which contain vertices of $T$.  The following conjecture is a natural generalization of Corollary \ref{nwt_cor} to $T$-Steiner trees.

\begin{conjecture}[Kriesell \cite{Kr1}]\label{kri}
Let $G$ be a connected graph and let $T \subseteq V(G)$.  If every edge-cut of $G$ separating $T$ has size at least $2k$ then
$G$ contains $k$ pairwise edge-disjoint $T$-Steiner trees.
\end{conjecture}

In contrast to spanning trees, $T$-Steiner trees are not the bases of any matroid, and hence the problem of packing them efficiently appears to be considerably more difficult.
Next we give a quick summary of partial results toward the above conjecture.  All of these results will be discussed in greater detail in the following section, where we describe the arguments used in proving our new theorem.

In the original article where his conjecture is posed, Kriesell shows that his conjecture is true in the special case when every vertex in $V \setminus T$ has even degree.  Frank, Kir\'aly and Kriesell \cite{FKK} have proved that there exist $k$ edge-disjoint $T$-Steiner trees whenever all edge cuts separating $T$ have size $\ge 3k$ and $V \setminus T$ is independent.  The first approximation to Kriesell's Conjecture without any added assumptions beyond edge-connectivity is due to Lau \cite{La} who proved that the conjecture holds true if the ``$2k$'' term is replaced by ``$24k$''.  More recently, West and Wu \cite{WW} have made a significant improvement on this, showing that the conjecture still holds when the ``$2k$'' is replaced by $``6.5k"$.  Our main theorem gives a further improvement on this bound.

\begin{theorem}
\label{steiner_thm}
Let $G$ be a graph and let $T \subseteq V(G)$. If every edge-cut in $G$ separating $T$ has size $\ge 5k+4$ then $G$ contains $k$ pairwise edge-disjoint $T$-Steiner trees.
\end{theorem}

For a pair of edges $e=uv$, and $e' = uv'$ we \emph{split} $e$ and $e'$ from $u$ by deleting $e$ and $e'$ and adding a new edge $vv'$.  West and Wu also introduced the following structure which is closely related to $T$-Steiner trees:  a $T$-\emph{connector} is a subgraph of $G$ which may be turned into a connected subgraph with vertex set $T$ by splitting pairs of edges. This requirement means that, for example, a $T$-connector cannot have a vertex of odd degree unless it is in $T$.  They make the following conjecture concerning the packing of $T$-connectors.

\begin{conjecture}[West and Wu \cite{WW}]
Let $G$ be a connected graph and let $T \subseteq V(G)$. If every edge-cut of $G$ separating $T$ has size $\ge 3k$ then $G$ contains $k$ pairwise edge-disjoint $T$-connectors.
\end{conjecture}

To see that this conjecture (if true) is essentially tight, consider the extreme circumstance where $G=(V,E)$ is a bipartite graph with bipartition $(T,V \setminus T)$ and every vertex in $V \setminus T$ has degree three.  If $|T| = n$, then each $T$-connector must contain at least $n-1$ vertices in $V \setminus T$ which do not appear in any other connector.  So in order for $G$ to contain $k$ edge-disjoint $T$-connectors we must have $|V \setminus T| \ge k(n-1)$.  In this case $|E| \ge 3k(n-1)$ and the average degree of a vertex in $T$ must be at least $3k \frac{n-1}{n}$.

West and Wu proved their conjecture holds if the ``$3k$'' term is replaced by $10k$.  We show the following improvement.

\begin{theorem}
\label{connector_thm}
Let $G$ be a graph and let $T \subseteq V(G)$. If every edge-cut in $G$ separating $T$ has size $\ge 6k+6$ then $G$ contains $k$ pairwise edge-disjoint $T$-connectors.
\end{theorem}

\section{Proof Discussion}

The proof of our main theorem is rather involved, but most of the essential ingredients can be traced back to previous works on Kriesell's Conjecture.  Accordingly, in this section we shall describe some of these earlier results, and indicate how we will utilize them.  The first key tool in our proof is the following famous theorem of Mader
(in fact Mader's theorem assumes the stronger condition that $u$ is not a cut vertex, but the theorem stated here follows from this).

\begin{theorem}[Mader's Splitting Theorem \cite{Ma}]
\label{mader_split}
Let $G$ be a graph, let $u \in V(G)$, and assume that ${\mathit deg}(u) \neq 3$ and that $u$ is not incident with a cut-edge.  Then there exists a pair of edges that can be split off from $u$ so that for all $x, x' \in V(G)\setminus u$, the size of the smallest edge-cut separating $x$ and $x'$ does not change.
\end{theorem}

In \cite{Kr1} Kriesell observes that Mader's Splitting Theorem can be used to verify his conjecture in the special case when every vertex in $V \setminus T$ has even degree.  Indeed, for such a graph, we may repeatedly apply Theorem \ref{mader_split} to each vertex in $V \setminus T$ until all such vertices have become isolated.  Deleting these vertices from the graph results in a $2k$-edge-connected graph with vertex set $T$ which has $k$ edge-disjoint spanning trees by Corollary \ref{nwt_cor}.  It is straightforward to verify that this yields $k$ edge-disjoint $T$-Steiner trees in the original graph.

Frank, Kir\'aly and Kriesell \cite{FKK} proved that Kriesell's Conjecture holds under the stronger assumptions that every edge-cut separating $T$ has size at least $3k$ and $V \setminus T$ is an independent set (in fact, they prove the stronger result that such a graph contains $k$ pairwise edge-disjoint $T$-connectors).  Since we will exploit these ideas further, let us sketch their argument.  Suppose that we are in this setting, and observe that we may repeatedly apply Mader's Splitting Theorem to reduce to the case when all vertices in $V \setminus T$ have degree $3$.  Since we may assume that a $T$-Steiner tree has no leaf vertex in $V \setminus T$, at this point, each vertex in $V \setminus T$ need only appear in one of the $T$-Steiner trees.  To model this, we now form a hypergraph $H$ with vertex set $T$ from our graph $G$ by replacing each $v \in V \setminus T$ by a single hyperedge on the three neighbours of $v$.  Now to solve the original problem, it suffices to find $k$ edge-disjoint spanning connected subhypergraphs of $H$.  Although spanning connected subhypergraphs of $H$ do not form a matroid, there is a matroid called the \emph{hypergraphic} matroid which is of use.  This is a matroid introduced by Lorea and defined on the edges of a hypergraph which we discuss in detail in the following section.  Assuming a hypergraph is suitably rich, any basis of the hypergraphic matroid will be a spanning connected subhypergraph.  Frank, Kir\'aly and Kriesell show that this hypergraph $H$ is suitably rich, and moreover it contains $k$ disjoint bases, and this gives them the desired Steiner Tree packing.

Lau found a way to reduce the general problem down to one in which $V \setminus T$ is independent, and thus take advantage of the hypergraphic matroid.  To motivate his approach, suppose that we have a graph $G = (V,E)$ with a set $T \subseteq V$ for which every edge-cut separating $T$ has size at least $f(k)$, in which we hope to prove the existence of $k$ edge-disjoint $T$-Steiner trees.  Let $e \in E$ and assume that $e$ is not incident with any vertex in $T$.  If $e$ is not in any edge-cut separating $T$ of size $f(k)$ then we may simply delete it.  So, let us assume that $\delta(X)$ is an edge-cut of $G$ separating $T$ with size $f(k)$ and $e \in \delta(X)$.  Since both $X$ and its complement contain a point in $T$ and an end of $e$, we have $|X|, |V \setminus X| \ge 2$.   So, a natural approach is to reduce this problem to smaller problems on the two graphs obtained from $G$ by identifying either $X$ or its complement to a single vertex which is placed in $T$ (note that each such graph still has the property that any edge cut separating its vertices in $T$ has size at least $f(k)$).

\begin{figure}[htb]
  \centering
  \includegraphics[width=200pt]{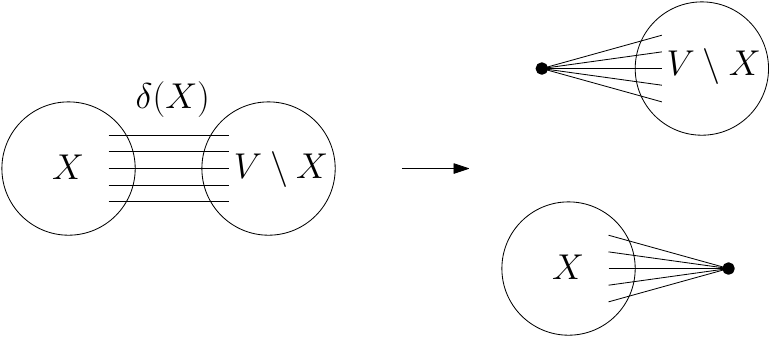}
  \caption{Reducing a nontrivial edge cut.}
  \label{figeasy}
\end{figure}

\begin{figure}[htb]
  \centering
  \includegraphics[width=200pt]{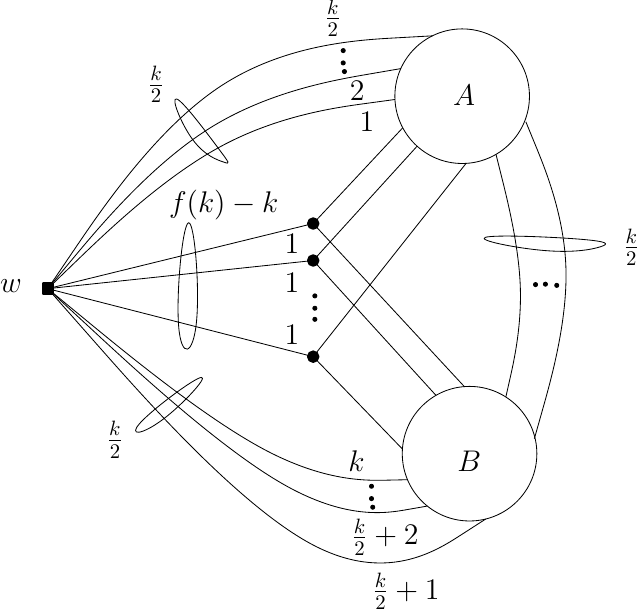}
  \caption{An obstruction to the na\"{\i}ve approach.}
  \label{figblanks}
\end{figure}

There are a number of difficulties with this approach.  First off we must arrange that the behaviour of our two solutions agree on $\delta(X)$; second,  we need to be sure that the resulting subgraphs we obtain are still connected.  A natural way to arrange this is to prove a stronger result where the behaviour of the packing is predetermined at a particular vertex.  So, let us now fix a distinguished vertex $w \in T$ of degree $f(k)$ and a function $\phi : \delta(w) \rightarrow \{0,1,\ldots,k\}$ and consider the problem of finding $k$ edge-disjoint connected subgraphs $G_1, \ldots, G_k$ with the following properties:
\begin{enumerate}
\item Every $G_i$ spans $T$
\item $E(G_i) \cap \delta(w) = \phi^{-1}(i)$ for every $1 \le i \le k$.
\item $G_i - w$ remains connected for every $1 \le i \le k$.
\end{enumerate}
These assumptions do allow us to naturally reduce any edge-cut $\delta(X)$ of size $f(k)$ with $|X|, |V \setminus X| \ge 2$ as suggested.  Unfortunately, there exist functions $\phi$ for which this problem has no solution, even when $f(k)$ is large.  One such case is indicated in Figure \ref{figblanks}.

In some sense, the difficulty in the example from Figure \ref{figblanks} is the fact that too many of the precoloured edges have the same label.  To get around this, Lau moves to yet another problem.  He considers the same problem as before with $G$ and $w$ and $\phi$, however he now attempts to pack $k$ edge-disjoint connected subgraphs $G_1, \ldots, G_k$ so that each spans $T$ and so that every vertex in $T$ has  degree at least two in every $G_i$.  This permits him the additional assumption that $|\phi^{-1}(i)| \ge 2$ for every $1 \le i \le k$, and Lau shows that for $f(k) = 26k$ this approach can be made to work.

Recently West and Wu extended these ideas further to prove that Kriesell's Conjecture holds true under the weaker assumption that all edge-cuts separating $T$ have size at least $6.5k$.  In fact, these authors prove a somewhat stronger statement which involves a more complicated simultaneous packing.  Their argument has some similarities to that of Lau in the sense that it features a special vertex at which the packing has been predetermined.  However, instead of insisting that each $G_i$ have minimum degree $\ge 2$, these authors try to pack the graphs $G_i$ so that every vertex in $T$ has degree at least $k$ in the graph $G \setminus \cup_{i=1}^k G_i$ (i.e. the packing leaves $k$ unused edges at each vertex in $T$).  This allows them to assume that the function $\phi$ satisfies $|\phi^{-1}(0)| \ge k$ and this gets around the obstruction in Figure \ref{figblanks}.

Our approach will follow a very similar line to that of West and Wu to reduce to a problem concerning finding many disjoint bases in a certain hypergraphic matroid.  The main novelty in our argument is that we pass to the fractional setting to find a packing, and then use the following pretty theorem to round this fractional solution to an integral one (which avoids using too many edges at any vertex in $T$)

\begin{theorem}[Kir\'aly, Lau, Singh \cite{KLS}]
\label{kls}
Let $M$ be a matroid on $E$, let $x \in {\mathbb R}^E$ be a fractional basis, and let ${\mathcal F}$ be a collection of subsets of $E$ so that every $e \in E$ is contained in at most $d$ members of ${\mathcal F}$.  Then there exists a basis $B$ so that every $F \in {\mathcal F}$ satisfies
\begin{equation}\label{klsintersect} |B \cap F| \le \lceil x(F) \rceil + d - 1.\end{equation}
\end{theorem}

The main advantage of passing to a fractional basis packing problem is that it permits us to take advantage of LP-duality (more precisely the Farkas' Lemma), and this
serves as the final ingredient in our proof.

\section{Matroid Preliminaries}

We shall require a couple of results concerning matroids.  The main matroid of interest to us will be the hypergraphic matroid associated with a hypergraph.  The essential properties of hypergraphic matroids we require appear already in Frank-Kir\'aly-Kriesell \cite{FKK}, but we have a slightly different approach toward the proofs which we feel has its own merit.

For every set $S$ we let $2^S$ denote the collection of all subsets of $S$.  A function $f : 2^S \rightarrow {\mathbb Z}$ is \emph{monotone increasing} if $f(X) \le f(Y)$ whenever $X \subseteq Y \subseteq S$ and is
\emph{intersecting submodular} if
\[ f(X \cap Y) + f(X \cup Y) \le f(X) + f(Y) \]
whenever $X,Y \subseteq S$ satisfy $X \cap Y \neq \emptyset$.  Edmonds and Rota~\cite{ER} proved the following fundamental result stating that every function $f$ with these two properties gives rise to a matroid on $S$.

\begin{theorem}[Edmonds and Rota \cite{ER}]
\label{int_sub}
Let $f : 2^S \rightarrow {\mathbb Z}$ be nonnegative with $f( \emptyset ) = 0$, and assume $f$ is monotone increasing and intersecting submodular.
Then there is a matroid $M$ on $S$ whose independent sets are
\[  \{ I \subseteq S \mid f(I') \ge |I'| \textrm{ for all } I'\subseteq I \} \]
and has rank function $r$ given by the rule
\[ r(A) = \min \Big\{ |R| +  \sum_{i=1}^{\ell} f(S_i) \mid \mbox{$\{R, S_1, \ldots, S_{\ell} \}$ is a partition of $A$} \Big\} .\]
\end{theorem}

Rado~\cite{Rado} proved a fundamental theorem which reveals how a bipartite graph with bipartition $(X,Y)$ together with a matroid on $Y$ induces the structure of a matroid on $X$.  His result may be deduced in quite general form from that of Edmonds above.  To do so, let $k,b \in {\mathbb Z}$ with $k \ge 1$ and define the set function $f : 2^X \rightarrow {\mathbb Z}$ by the rule $f(S) = \max\{ 0,  k r (N(S)) + b \}$.  Then $f$ is monotone increasing and intersecting submodular, so it gives rise to a matroid on $X$ by Theorem \ref{int_sub} whose nonempty independent sets are all $\emptyset \neq I \subseteq X$ which satisfy
\[  |I'| \le k r(N(I')) + b, \textrm{\qquad for all } I'\subseteq I. \]
We call this matroid the \emph{Rado} matroid on $X$ (though the reader should note that every matroid $L$ may be obtained in this manner by taking both $X$ and $Y$ to be copies of $E(L)$ and letting the graph be the matching where $x \sim y$ if $x$ and $y$ are copies of the same element in $E(L)$).  The \emph{standard} Rado matroid has $k=1$ and $b=0$ and these matroids are especially nice.
Notably, transversal matroids are the standard Rado matroids where the matroid on $Y$ is free (i.e. all sets are independent), and the bicircular matroid of a graph $G=(V,E)$ is the transversal matroid coming from the associated bipartite graph with bipartition $(E,V)$ (where $ev$ is an edge if $e$ is incident with $v$).  Using this same bipartite graph with bipartition $(E,V)$ associated with the graph $G$ but taking $b = -1$ gives us the usual cycle matroid of the graph $G$.  More generally, if $H = (V,F)$ is a hypergraph and we take the bipartite graph with bipartition $(F,V)$ (where $fv$ is an edge if $v \in f$) the free matroid on $V$ and the parameters $k=1$ and $b=-1$, we get the \emph{hypergraphic} matroid of the hypergraph $H$, first discovered by Lorea~\cite{Lorea}.  However, this matroid has an alternate representation as a standard Rado matroid.  Namely, we define $E$ to be the edge set of a complete graph on $V$, define the bipartite graph with bipartition $(F,E)$ by the rule that each $e \in F$ is adjacent to every edge of the form $uv$ where $u,v$ are contained in $e$, and use the parameters $k=1$ and $b=0$ with the standard cycle matroid of the graph $(V,E)$ as the matroid on $E$.  Indeed, this representation may be viewed as a halfway step between $F$ and $V$ by constructing a tripartite graph on $F \cup E \cup V$ with the natural adjacencies.

\bigskip

The alternate representation has some advantages for us.  Perhaps the main advantage is that it gives us a representation of the hypergraphic matroid of $H$ as a standard Rado matroid.  For standard Rado matroids, Rado proved that a subset of $X$ is independent if and only if it can be matched to an independent set in $Y$ and this immediately yields the following result.

\begin{theorem}[Frank, Kir\'aly, Kriesell \cite{FKK}]\label{thm22}
If $H = (V,F)$ is a hypergraph and $A \subseteq F$ then $A$ is independent in the hypergraphic matroid if and only if we can
replace every $e \in A$ with an edge $uv$ where $u,v \in e$ so that the resulting graph is a forest.
\end{theorem}

Next let us turn our attention to the rank function of a hypergraphic matroid on $H = (V,E)$.  Here it is convenient to work in the language
of partitions.  For every partition ${\mathcal P}$ of $V$ we define the \emph{rank} of ${\mathcal P}$ to be $r( {\mathcal P} ) = |V| - |{\mathcal P}|$.
The following result gives a nice form for the rank function of the hypergraphic matroid.  This can be proved rather easily using either of the two Rado representations of this matroid.  To deduce it from the representation as a standard Rado matroid on $(F,E)$, note that for every $A \subseteq F$
we have $r(N(A))$ is equal to $r( {\mathcal P} )$ where ${\mathcal P}$ is the partition of $V$ induced by the connected components of $A$.  Since we have $\sum_{i=1}^{\ell} r(N(S_i)) \ge r ( N( \cup_{i=1}^{\ell} S_i ))$, the next result follows from the rank formula in Theorem \ref{int_sub}.

\begin{theorem}[Frank, Kir\'aly, Kriesell \cite{FKK}]\label{steiner_basis}
In the hypergraphic matroid on a hypergraph $H = (V,F)$ we have for every set $A \subseteq F$ that the rank of $A$ is equal to the
minimum over all partitions ${\mathcal P}$ of $V$ of the function
\[ r({\mathcal P}) + \lambda^{out}_{\mathcal P}(A) .\]
\end{theorem}

Another convenience of the present development is that matroid union is very easy to deal with as indicated by the following observation.
For every matroid $M$ we let $kM$ denote the matroid union of $k$ copies of $M$.

\begin{observation}
Let $L$ be a Rado matroid associated with the bipartite graph on $(X,Y)$, the matroid $M$ on $Y$, and the parameters $m,b$. Then $kL$ is the same Rado matroid as $L$ except with the parameters $km$ and $kb$ in place of $m$ and $b$.
\end{observation}

\noindent{\it Proof:} Using the formula for matroid union yields that for every $A \subseteq E(M)$
\begin{align*}
r_{kM}(A)
	&= \min_{A_1 \subseteq A} \Big(|A \setminus A_1| + k r_M(A_1)\Big)	\\
	&= \min_{A_1 \subseteq A} \Big(|A \setminus A_1| + k \min_{A_2 \subseteq A_1}  \Big( |A_1 \setminus A_2| + \max\{0, m r(N(A_2)) + b \} \Big)\Big) 	\\
	&= \min_{A^* \subseteq A} \Big(|A \setminus A^*| + k \max\{0, m r (N(A^*)) + b \}\Big)	\\
	&= \min_{A^* \subseteq A} \Big(|A \setminus A^*| + \max\{0, km r (N(A^*)) + kb \}\Big)
\end{align*}
\vspace*{-.5in}

\hfill$\Box$

\bigskip

Accordingly, we obtain the following formula for the rank function in a union of $k$ copies of a hypergraphic matroid.

\begin{theorem}[Frank, Kir\'aly, Kriesell \cite{FKK}]
\label{union_rank}
In the union of $k$ copies of the hypergraphic matroid on a hypergraph $H = (V,F)$, we have for every set $A \subseteq F$ that the rank of $A$ is equal to the
minimum over all partitions ${\mathcal P}$ of $V$ of the function
\[ k r({\mathcal P}) + \lambda^{out}_{\mathcal P}(A). \]
\end{theorem}

There are two important polyhedra associated with a matroid and we shall require some basic properties of these as well.  For any matroid $M$ on $E$ consider the following three linear constraints on vectors $x \in {\mathbb R}^E$.
\begin{enumerate}
\item[(L1)] $0 \le x(e) \le 1$ for every $e \in E$.
\item[(L2)] $x(S) \le r(S)$ for every $S \subseteq E$.
\item[(L3)] $x(E) = r(E)$.
\end{enumerate}
Edmonds~\cite{Edmonds} proved that a vector $x \in {\mathbb R}^E$ satisfies (L1) and (L2) above if and only if $x$ can be expressed as a convex combination of incidence vectors of independent sets, and we call this polyhedron ${\mathcal I}(M)$ the \emph{independent set polytope} of $M$.  It follows from this that $x$ satisfies (L1), (L2), and (L3) if and only if $x$ can be expressed as a convex combination of incidence vectors of bases, and we call this polyhedron
${\mathcal B}(M)$ the \emph{basis polytope} of $M$.  We call an element of ${\mathcal B}(M)$ a \emph{fractional basis} of $M$.  The key property we will require in the following sections is the following.

\begin{lemma}
\label{partition_rank}
Let $M$ be a hypergraphic matroid associated with a hypergraph on $V$ and suppose that $x \in {\mathbb R}^E$ satisfies the following
\begin{enumerate}
\item $0 \le x(e) \le 1$ for every $e \in E$.
\item $x( \Lambda^{in}_{\mathcal P}(E ) ) \le k r( {\mathcal P} )$ for every partition ${\mathcal P}$ of $V$.
\item $x(E) = k( |V| - 1)$.
\end{enumerate}
Then $kM$ has rank $k r(M)$ and $x$ is a fractional basis of $kM$.
\end{lemma}

\noindent{\it Proof:} Let $A \subseteq E$ and choose a partition ${\mathcal P}$ in accordance with Theorem \ref{union_rank}.
This gives us
\begin{align*}
r(A)	& =	k r(  {\mathcal P} ) + \lambda^{out}_{\mathcal P}(A)	\\
	&\ge x( \Lambda^{in}_{\mathcal P}(A) ) + x( \Lambda^{out}_{\mathcal P}(A) )	\\
	& = x(A)
\end{align*}
Therefore, $x$ may be expressed as a convex combination of incidence vectors of independent sets.  Since $x(E) = k (|V| - 1)$  and $r(M) \le |V|-1$
it follows that $kM$ must have rank $k r(M) = k (|V| - 1)$ and $x$ is a fractional basis of $kM$ as desired.\hfill$\Box$

\bigskip

We require just one additional result which will permit us to arrange a small number of edges as desired within a matroid union.

\begin{lemma}
\label{adjust_union}
Let $M$ be a matroid, let $B$ be independent in $kM$, let $D \subseteq B$ and let $\mu : D \rightarrow \{1,\ldots,k\}$ be an injection.  Then there exist pairwise disjoint sets $B_1, \ldots, B_k$ satisfying
\begin{itemize}
\item $\cup_{i=1}^k B_i = B$.
\item Every $B_i$ is independent in $M$.
\item $\mu^{-1}( \{i\} ) \subseteq B_i$ for $1 \le i \le k$.
\end{itemize}
\end{lemma}

\noindent{\it Proof:}
Since $B$ is independent in $kM$ we may choose pairwise disjoint sets $B_1,\ldots,B_k \subseteq B$ so that each $B_i$ is independent in $M$ and $\cup_{i=1}^k B_i = B$.  We shall now show how to modify these sets so that they satisfy the last condition.  If there exists $1 \le i \le k$ with
$\mu^{-1}( \{i\} )$ not contained in $B_i$ then we must have $\mu^{-1}( \{i\} ) = \{e\}$ and $e \in B_j$ for some $j \neq i$.  If $B_i \cup \{e\}$ is independent in $M$ then we make this set our new $B_i$ and set $B_j = B_j \setminus \{e\}$.  Otherwise, we may choose $f \in B_i$ so that both $(B_i \cup \{e\}) \setminus \{f\}$ and $(B_j \cup \{f\}) \setminus \{e\}$ are independent in $M$ and we define $B_i$ and $B_j$ respectively to be these new independent sets.  It follows easily that this modification increases the number of indices $1 \le i \le k$ with $\mu^{-1}( \{i\} ) \subseteq B_i$ so it will terminate with a solution.\hfill$\Box$

\section{Packing Steiner Trees}\label{sec:steiner}

The goal of this section is to establish Theorem \ref{steiner_thm}.  In order to do so, we will prove the following
stronger result which features a distinguished vertex whose incident edges have been precoloured.
This lemma involves the function $f : {\mathbb Z} \rightarrow {\mathbb Z}$ which we define as follows:
\[ f(k) = 2 \left\lceil \frac{5k+3}{2} \right\rceil. \]

\begin{lemma}
\label{steiner_lem}
Let $k \ge 3$, let $G$ be a graph, and let $T \subseteq V(G)$ satisfy $|T| \ge 2$ and have the property that every edge cut separating $T$ has size
$\ge f(k)$.  Let $w \in T$ be a distinguished vertex with ${\mathit deg}(w) = f(k)$ and let $\phi : \delta(w) \rightarrow \{0,1,\ldots,k\}$
be onto and satisfy $| \phi^{-1}(0) | \ge k$.  Then there exist $k$ edge-disjoint connected subgraphs $G_1,\ldots,G_k$ that satisfy
all of the following conditions:
\begin{enumerate}
\item[(C1)] $G_i$ spans $T$ for $1 \le i \le k$.
\item[(C2)] $\delta(w) \cap E(G_i) = \phi^{-1}(i)$ for $1 \le i \le k$.
\item[(C3)] $w$ is not a cut vertex of $G_i$ for $1 \le i \le k$.
\item[(C4)] Every $v \in T$ satisfies $\sum_{i=1}^k {\mathit deg}_{G_i}(v) \le {\mathit deg}_G(v) - k$.
\end{enumerate}
\end{lemma}

Lemma \ref{steiner_lem} immediately implies Theorem \ref{steiner_thm}, as follows.

\bigskip

\noindent{\textbf{\it Proof of Theorem \ref{steiner_thm}:}} We may assume that $k \ge 3$ as otherwise the result follows from the earlier result of West and Wu \cite{WW}.  For this case we shall prove the stronger result that the same conclusion holds under the assumption that all edge-cuts separating $T$ have size at least $f(k)$.  By possibly deleting edges, we may assume that there exists an edge-cut $D$ in $G$ which separates $T$ and has size $f(k)$.  Now, choose a function $\phi : D \rightarrow \{0,1,\ldots,k\}$ which is onto and satisfies $| \phi^{-1}(0)| \ge k$.  Let $\{X_1, X_2\}$ be the vertex partition associated with $D$ and form two new graphs as follows:  form $G^1$ from $G$ by identifying $X_2$ to a new vertex called $x_2$ and form $G^2$ from $G$ by identifying $X_1$ to a new vertex called $x_1$.  Now apply Lemma \ref{steiner_lem} to the graph $G^1$ with $T^1 := (T \cap X_1) \cup \{x_2\}$, the distinguished vertex $x_2$ and the function $\phi$, to get $k$-edge-disjoint subgraphs $G^1_1,\ldots, G^1_k$ of $G^1$ satisfying (C1)-(C4) for the given parameters. Similarly, apply Lemma \ref{steiner_lem} to $G^2$ for $T^2=(T \cap X_2) \cup \{x_1\}$, the vertex $x_1$ and the function $\phi$ to get $k$-edge-disjoint subgraphs $G^2_1,\ldots, G^2_k$ of $G^2$ satisfying (C1)-(C4) for the given parameters. Combine $G^1_i$ and $G^2_i$ in the obvious way to get a subgraph $G_i$ of $G$, for all $i\in \{1,2,\ldots,k\}$. By (C1) and (C3), and since $\phi$ is onto, we get that each such $G_i$ is connected and spans $T$. Hence each such $G_i$ contains a $T$-Steiner tree.
\hfill$\Box$

\bigskip

Note that for the purposes of the above proof, outcome (C4) and the assumption of $\phi^{-1}(0)\geq k$ in Lemma \ref{steiner_lem} are extraneous. This idea of leaving at least $k$ edges at each vertex of $T$ which are not used in the Steiner trees will however be key in establishing Lemma \ref{steiner_lem}. Our proof of Lemma \ref{steiner_lem} will involve investigating a minimal counterexample, and we shall instantiate this here.

For the remainder of this section, assume that $G$ together with $T \subseteq V(G)$, $w \in T$ and $\phi$ is a
counterexample to Lemma \ref{steiner_lem} which is extreme in the following sense:
\begin{enumerate}
\item[(X1)] $\sum_{v \in T} {\mathit deg}(v)$ is minimum.
\item[(X2)] $|V(G)| + |E(G)|$ is minimum (subject to X1).
\item[(X3)] $|\phi^{-1}(0)|$ is minimum (subject to X2).
\end{enumerate}

Our first goal will be to establish the following regularity and independence properties of the graph $G$.

\begin{lemma}
\label{minimal_counter}
The graph $G$ satisfies the following properties:
\begin{itemize}
\item[(P1)] $|\phi^{-1}(0)| = k$.
\item[(P2)] Every $v \in V(G) \setminus T$ has degree 3 and has $|N(v)| = 3$.
\item[(P3)] Every $v \in T$ has degree $f(k)$.
\item[(P4)] $V(G) \setminus T$ is an independent set.
\item[(P5)] $|T| \ge 3$.
\end{itemize}
\end{lemma}

\noindent{\it Proof:} We shall establish the desired properties in a series of claims.

\bigskip

\noindent{0.} $|\phi^{-1}(0)| = k$.

\smallskip

By assumption we must have $|\phi^{-1}(0)| \ge k$.  If this set were larger, then we could obtain a smaller counterexample (via (X3)) by choosing $e \in \delta(w)$ with $\phi(e) =0$ and instead setting $\phi(e) = 1$.
(If the theorem holds after this adjustment, then by (C3) the graph $G_1$ will still be connected after removing the edge $e$.)

\bigskip

\noindent{1.} $G$ is connected and has no cut-edge.

\smallskip
Suppose (for a contradiction) that $G$ is disconnected.  Then one component of $G$ must contain all of $T$, and this component forms a smaller counterexample (via (X2)), giving us a contradiction.
Now suppose (for a contradiction) that $e \in E(G)$ is a cut-edge.  Then one component of $G-e$ must contain all of $T$, and it follows that $e$ cannot be in any edge-cut of size $f(k)$ which separates $T$.  However, then the component of $G-e$ containing $T$ forms a smaller counterexample (via (X1) or (X2)), giving us a contradiction.

\bigskip

\noindent{2.} Every $v \in V(G) \setminus T$ has degree 3.

\smallskip

Let $v \in V(G) \setminus T$ and note that by (1) we must have ${\mathit deg}(v) \ge 2$.  If ${\mathit deg}(v) = 2$, then splitting off the two edges at $v$ and suppressing $v$ gives a smaller counterexample (via (X2)), so we must have ${\mathit deg}(v) \ge 3$.  If ${\mathit deg}(v) \ge 4$ then by applying Mader's Splitting Theorem we obtain a new graph $G'$ which will again be a smaller counterexample (via (X2)).  It follows that ${\mathit deg}(v) = 3$ as desired.

\bigskip

\noindent{3.} Every $v \in V(G) \setminus T$ has $|N(v)| = 3$.

\smallskip

If this property were false, then by (2) there would be a (parallel) edge $e$ incident with $v$ which is not contained in any edge-cut separating $T$ of size $f(k)$.  However, then $G-e$ is a smaller counterexample (via (X1) or (X2)).

\bigskip

\noindent{4.} If $D$ is an edge-cut separating $T$ and $|D| = f(k)$ then $D = \delta(v)$ for some $v \in T$.

\smallskip

Suppose (for a contradiction) that $D$ partitions $V(G)$ into two sets $X_1,X_2$ with $|X_i| \ge 2$ for $i=1,2$.  We may assume without loss of generality that $w \in X_1$. We now form two new graphs, $G^1$ and $G^2$, exactly as described in the proof of Theorem \ref{steiner_thm} above, and also  define $T^1$ and $T^2$ as in this proof. We now apply induction to $G^1$, $T^1$, $\phi$ and distinguished vertex $w$; our induction is via (X1) or (X2). This application yields subgraphs $G^1_1,\ldots,G^1_k$ of $G^1$ satisfying (C1)-(C4) for the given parameters.

Now  let the map $\phi_2 : D \rightarrow \{0,1,\ldots,k\}$ be defined by the rule that $\phi_2^{-1}(i) = D \cap E(G^1_i)$ for $1 \le i \le k$, and $\phi_2^{-1}(0)=D\setminus \cup_{i=1}^k E(G^1_i)$. Since $x_2 \in T^1$ it must be that $\phi_2$ is onto and further (by (C4)) that $|\phi_2^{-1}(0)| \ge k$.
Now we apply induction to $G^2$, $T^2$, the distinguished vertex $x_1$, and the map $\phi_2$; our induction is via (X1) or (X2). This application yields subgraphs $G^2_1,\ldots,G^2_k$ of $G^2$ satisfying (C1)-(C4) for the given parameters. Combine $G^1_i$ and $G^2_i$ in the obvious way to get a subgraph $G_i$ of $G$, for all $i\in \{1,2,\ldots,k\}$. These subgraphs now satisfy (C1)-(C4) for the original parameters of our problem, contradicting the fact that $G$ is a counterexample to Lemma \ref{steiner_lem}.

\bigskip

\noindent{5.} Every $v \in T$ has degree $f(k)$.

\smallskip

Suppose (for a contradiction) that $v \in T$ satisfies ${\mathit deg}(v) > f(k)$ and note that $v \neq w$ by assumption.  Now apply Mader's Splitting Theorem to choose two edges $uv$ and $u'v$ incident with $v$ so that the graph $G'$ obtained from $G$ by splitting these edges at $v$ preserves edge-connectivity between pairs of vertices not including $v$.  Then, we let $G''$ be the graph obtained from $G - \{uv,u'v\}$ by adding a new vertex $x$ and adding the edges $vx,ux,u'x$. We will show that every edge-cut in $G''$ separating $T$ has size $\geq f(k)$, which  implies that $G''$ (together with $T, w, \phi$) is a smaller counterexample (via (X1)), thus giving a contradiction. To this end, first note that every edge-cut of $G'$ which separates $T \setminus \{v\}$ has size $\ge f(k)$, so this property also holds for the graph $G''$. By (3) and the fact that $G$ is a counterexample we know that $|T|\geq 3$. By (4) we know that in the original graph $G$, every edge-cut which separates $\{v\}$ from $T \setminus \{v\}$ must have size $> f(k)$. Hence, in $G'$, every edge-cut which separates $\{v\}$ from $T \setminus \{v\}$ will have size $\ge f(k)$, and this property also holds for $G''$. This yields the desired contradiction.

\bigskip

\noindent{6.} $V(G) \setminus T$ is an independent set.

\smallskip

If $e$ is an edge with both ends in $V(G) \setminus T$ then it follows from (4) that $e$ is not in any edge-cut separating $T$
with size $f(k)$.  But then the graph $G-e$ (together with $T$, $w$, and $\phi$) is a smaller counterexample (via (X2)).

\bigskip

\noindent{7.} $|T| \ge 3$

\smallskip

If this were to fail we must have $|T| = 2$.  However, then (3) and (6) imply that $V(G) = T$, and it is easy to check that there is no counterexample on 2 vertices.  \hfill$\Box$

\bigskip

The function $\phi$ has specified the behaviour of any solution at the distinguished vertex $w$, and we next introduce some notation to help us to deal with the many possibilities for this map.
 Define $W^{\circ}_t$ to be the set of edges in $\phi^{-1}(0)$ which have both ends in $T$ and let $W^{\circ}_n = \phi^{-1}(0) \setminus W^{\circ}_t$.  Since $\phi$ is onto, we may choose disjoint subsets $W^*_t$ and $W^*_n$ of $\delta(w) \setminus \phi^{-1}(0)$ so that $\phi$ maps $W^*_t \cup W^*_n$ bijectively to  $\{1,\ldots,k\}$ and so that every edge in $W^*_t$ has both ends in $T$ and every edge in $W^*_n$ has one end in $V(G) \setminus T$.  Let $W^{\bullet}_t$ denote the set of edges incident with $w$ and another vertex in $T$ which are not contained in $W^{\circ}_t \cup W^*_t$ and similarly, let $W^{\bullet}_n$ denote the set of edges incident with $w$ and with a vertex in $V(G) \setminus T$ which are not contained in $W^{\circ}_n \cup W^*_n$.

\begin{figure}
  \centering
  \includegraphics[width=380pt]{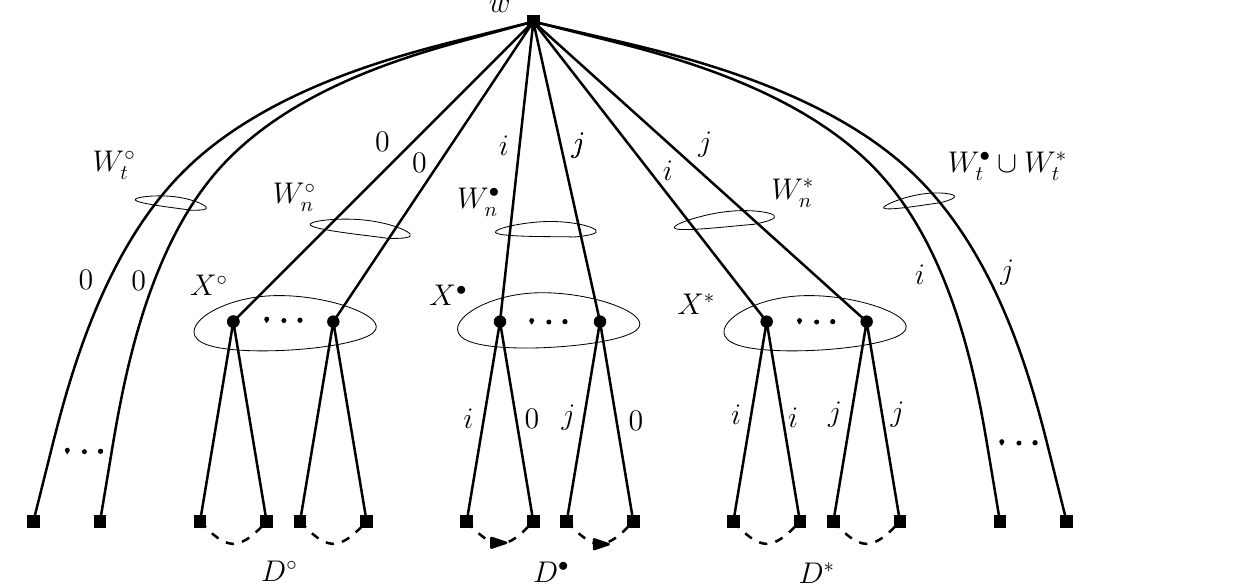}
  \caption{The distinguished vertex $w$ and its incident edges.}
  \label{figSpecialVertex}
\end{figure}

For each choice of $\star \in \{ \circ, *, \bullet \}$ let $X^{\star}$ be the set of endpoints of edges of the form $W_n^{\star}$ other than $w$.  So all vertices in $X^{\circ}$, $X^*$, and $X^{\bullet}$ are contained in $V(G) \setminus T$.  Next  we construct three sets of auxiliary edges $D^{\circ}$, $D^*$, and $D^{\bullet}$ which will play somewhat different roles in our proof.  For each choice of $\star \in \{ \circ, *, \bullet \}$ we construct $D^{\star}$ from $\emptyset$ as follows.  For each edge $x \in X^{\star}$ let $v,v'$ denote the neighbours of $x$ other than $w$ and add a new edge with ends $v$ and $v'$ to $D^{\star}$.  Note that this construction yields a natural correspondence between $X^{\star}$ and $D^{\star}$. See Figure \ref{figSpecialVertex}. Note that in this figure, and in later figures, we use squares to denote vertices in $T$, and circles to denote vertices not in T.

\bigskip

Our next step will be to define an extension $\phi'$ of our precolouring $\phi$, by defining $\phi'(e)=\phi(e)$ for all $e\in \delta(w)$, but also defining $\phi'$ on the other edges  of $\delta(x)$ where $x \in X^* \cup X^{\bullet}$.  First, for every $x \in X^*$ with $N(x) = \{w,v,v'\}$ we define $\phi'(xv) = \phi'(xv') = \phi(wx)$.  In order to extend our precolouring in the case when $x \in X^{\bullet}$ we need to be a little careful.  Define $G^{\bullet}$ to be the graph $(V(G), D^{\bullet})$ and then choose an orientation of $G^{\bullet}$ so that every vertex $v$ satisfies
  \begin{equation}
  \label{orientation1}
 {\mathit deg}^+_{G^{\bullet}}(v) \le {\mathit deg}^-_{G^{\bullet}}(v) + 1.
 \end{equation}
 (such an orientation may be found greedily as follows: repeatedly choose a cycle of unoriented edges and orient them cyclicly until the remaining unoriented edges form a forest, and then repeatedly orienting leaf edges of the forest of unoriented edges away from the leaf vertex).  Now, for each $x \in X^{\bullet}$ with $N(x) = \{w,v,v'\}$ we consider the edge $vv' \in D^{\bullet}$.  If this edge is oriented from $v$ to $v'$ then we define
\begin{align*}
	\phi'(xv) &= \phi(wx)	\\
	\phi'(xv') &= 0
\end{align*}

In proving Lemma \ref{steiner_lem}, we will reduce our original problem to a new one on a smaller graph which is obtained from the original by removing the vertex $w$ and edges precoloured by $\phi'$ and making some further modifications.  In preparation for this, we now introduce three functions which will keep track of the number of edges in the domain of $\phi'$ which are incident with a given vertex. For every $v \in T' = T \setminus \{w\}$ we define the following
\begin{align*}
	p(v)		&= | \delta(v) \cap (W^*_t \cup W^{\bullet}_t)| + | \{ x \in X^{*} \mid \mbox{$x$ is adjacent to $v$} \} |	\\
	q(v)		&= | \delta(v) \cap\phi'^{-1}(\{1,\ldots,k\}) | - p(v)	\\
	u(v)		&= | \delta(v) \cap \phi'^{-1}( \{0\} ) |
\end{align*}
Here $u(v)$ counts the number of edges incident with $v$ which are assigned the colour $0$ by $\phi'$ (the other end of such an edge could be either $w$ or a vertex in $X^{\bullet}$).  The quantity $p(v) + q(v)$ counts the number of edges incident with $v$ which are assigned a nonzero colour by $\phi'$.  These edges have three types, their endpoint other than $v$ might be $w$, might be in $X^*$ or might be in $X^{\bullet}$.  Here $p(v)$ counts the edges of the first two types, and $q(v)$ counts the edges of the third type.  Next we give an easy observation and a somewhat complicated lemma which give key bounds on $p(v)$ and $q(v)$.

\begin{observation}\label{obp2u} Every $v \in T'$ satisfies
\begin{itemize}
\item $q(v) \le u(v) + 1$
\item $2 q(v) \le f(k) + 1$
\end{itemize}
\end{observation}

\noindent{\it Proof:} The first inequality follows from the fact that the chosen orientation of $G^{\bullet}$ satisfies (\ref{orientation1}) and our subsequent definition of $\phi'$.   For the second we use the first to obtain $f(k) = \mathit{deg}(v) \ge q(v) + u(v) \ge 2q(v) - 1$.   \hfill$\Box$

\begin{lemma}
\label{buzzer_beater}
Every $v \in T'$ satisfies $p(v) \le k-1$.
\end{lemma}

\noindent{\it Proof:} Suppose (for a contradiction) that $p(v) \ge k$ and define the following set of
edges and set of colours.
\begin{align*}
A	&=	\{ e \in \delta(v) \mid \mbox{ $e$ is incident with a vertex in $X^*$} \}	\\
J	&=	\{ \phi'(e) \mid e \in A \}
\end{align*}
Note that our definitions imply $|J| = |A|$.  Since $p(v) \ge k$ we may choose a set $B \subseteq \delta(v) \cap (W_t^* \cup W_t^{\bullet})$ so that $|A| + |B| = k$ and we define $C = (\delta(v) \cap \delta(w)) \setminus B$.


Let $\hat{G}$ be the graph obtained from $G$ by applying Mader's Splitting Theorem repeatedly to split the vertex $v$ completely (this is possible since $deg_G(v)=f(k)$ by (P3) in Lemma~\ref{minimal_counter}, and $f(k)$ is always even).  If two edges in $B \cup C$ were split off together, the resulting graph would have a loop at $w$ and would then have an edge cut separating $w$ from $T \setminus \{v,w\}$ (which is nonempty by (P5)) of size smaller than $f(k)$ (contradicting Mader).  Similarly, if an edge in $B \cup C$ was split off with an edge in $A$, then in the resulting graph, some vertex $x \in X^*$ would have two edges between it and $w$, and the edge cut separating $w$ and $x$ from the remaining vertices would have size $< f(k)$ contradicting Mader.  Therefore, no such splits occur.


Let $\hat{B}$  and $\hat{C}$  denote  the set of edges in $\hat{G}$ which were formed by a split at $v$ involving an edge in $B$ and $C$ respectively, and note that $|\hat{B}|=|B|$ and $|\hat{C}| = |C|$.
We now construct a function $\hat{\phi} : \delta_{\hat{G}}(w) \rightarrow \{0,1,\ldots,k\}$ as follows.
We let $\hat{\phi}(e) = \phi(e)$ for every edge $e \in \delta_{\hat{G}}(w) \setminus (\hat{B} \cup \hat{C})$, we let $\hat{\phi}(e) = 0$ for every $e \in \hat{C}$ and we choose values for $\hat{\phi}(e)$ for those $e \in \hat{B}$ so that $\hat{\phi}(\hat{B}) = \{1,\ldots,k\} \setminus J$.

It now follows from the minimality of our counterexample, that Lemma \ref{steiner_lem} holds for the graph $\hat{G}$, the vertex subset $T \setminus \{v\}$, the distinguished vertex $w$ and the function $\hat{\phi}$.  This yields  a family of edge-disjoint connected subgraphs $\hat{G}_1, \ldots, \hat{G}_k$ of $\hat{G}$ satisfying (C1)-(C4) for the given parameters.  We will assume that these graphs are minimal in the sense that removing any edge from any $\hat{G}_i$ would violate one of the conditions.

Next we will start from $\hat{G}$ and return to the original graph $G$ by reversing our splitting operations one at a time, and modifying the graphs $\hat{G}_i$ along the way to obtain a solution for the original graph.  We begin by adding $v$ back to the graph $\hat{G}$ as an isolated vertex.

\begin{figure}
  \centering
  \includegraphics[width=250pt]{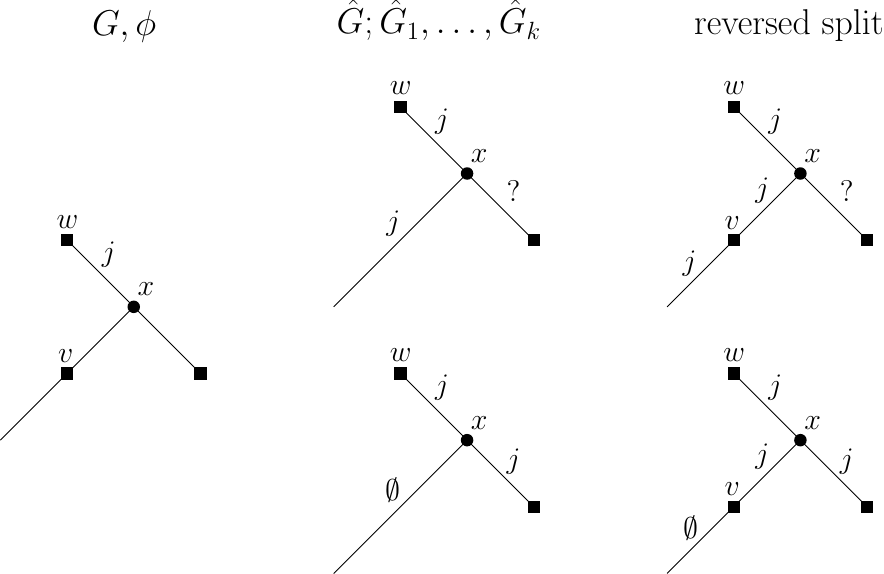}
  \caption{reversing the split of an $A$ edge (case 1).}
  \label{bbfig1}
\end{figure}


First consider an edge $e \in A$ and note that this was split off with some other edge $e' \in \delta(v) \setminus (B \cup C)$ to form an edge $\hat{e}$ in the graph $\hat{G}$.  In this case $e$ is incident with a vertex $x \in X^*$ and let us first consider the case that $e'$ is not also incident with a vertex of $X^*$.  By assumption, the edge between $x$ and $w$ in $\hat{G}$ is mapped by $\hat{\phi}$ to a (nonzero) colour $j \in J$.  Since $x$ has degree 3 and the graph $\hat{G}_j$ is minimal, we may assume that the graphs $\hat{G}_i$ with $i \neq j$ do not contain the vertex $x$.  If $\hat{e} \in E(G_j)$ then we delete $\hat{e}$ from $\hat{G}$ and add the edges $e,e'$ back to $\hat{G}$, and similarly we delete $\hat{e}$ from $\hat{G}_j$ and add the edges $e,e'$ to $\hat{G}_j$.  On the other hand, if $\hat{e} \not\in E(\hat{G}_j)$ then the edge incident with $x$ which is not $\hat{e}$ and not incident with $w$ must appear in the graph $\hat{G}_j$ (since $x \in V( \hat{G}_j)$ and $w$ is not a cut vertex of $\hat{G}_j$).  In this case we delete $\hat{e}$ and add the edges $e,e'$ back to our graph $\hat{G}$ and we add the edge $e$ to $\hat{G}_j$.  Note that after this operation, the graph $\hat{G}_j$ will span $T$ and will not have $w$ as a cut vertex.  This process appears in Figure \ref{bbfig1}; here the edge marked $j$ in the leftmost graph is an edge which is mapped to $j$ by the precolouring $\phi$, the edges marked $j$ in the middle column are those which appear in the graph $\hat{G}_j$, and the edges marked $j$ in the rightmost column are still in the graph $\hat{G}_j$ after we reverse the split.

\begin{figure}[h]
  \centering
  \includegraphics[width=250pt]{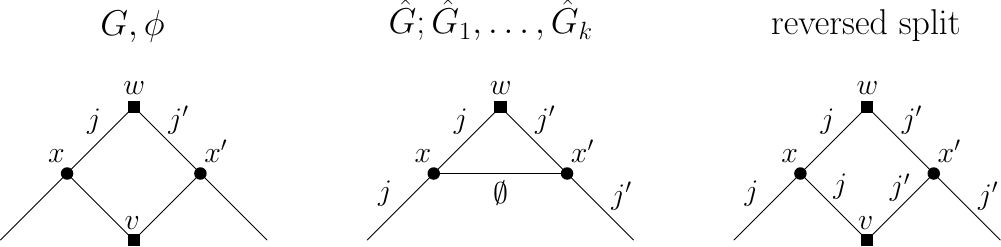}
  \caption{reversing the split of an $A$ edge (case 2).}
  \label{bbfig2}
\end{figure}

The other case is when the edge $e'$ is also incident with a vertex $x' \in X^*$ (see Figure \ref{bbfig2}).  Define  $j' \in J$ to be the value of $\hat{\phi}$ evaluated at the edge $wx'$.  Since $j \neq j'$ (by definition of $X^*$), our minimality assumptions imply that the edge $\hat{e}$ does not appear in any of the graphs $\hat{G}_i$.  However, by assumption we will have $x \in \hat{G}_j$ and $x' \in \hat{G}_{j'}$.  It now follows (from (C1) and (C3)) that the edge of $\hat{G}$ incident with $x$ $(x')$ which is not $\hat{e}$ and not incident with $w$ must appear in $\hat{G}_j$ $(\hat{G}_{j'})$.  Now we move back toward our original graph by deleting $\hat{e}$ and adding back the edges $e,e'$ and we add $e$ to $\hat{G}_j$ and $e'$ to $\hat{G}_{j'}$.  After this operation, the graphs $\hat{G}_j$ and $\hat{G}_{j'}$ will both span $T$ and not have $w$ as a cut vertex.

\begin{figure}[h]
  \centering
  \includegraphics[width=200pt]{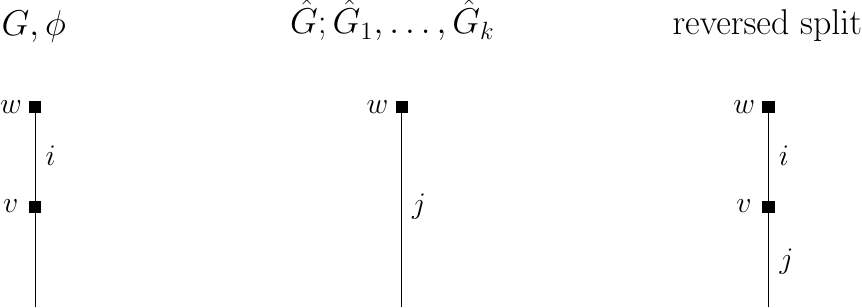}
  \caption{reversing the split of a $B$ edge.}
  \label{bbfig3}
\end{figure}

Next we will consider edges in $B$.  Let  $e \in B$, and suppose that $e$ was split off with the edge $e'$ to form the edge $\hat{e} \in \hat{B}$ (see Figure \ref{bbfig3}).  We modify $\hat{G}$ by deleting the edge $\hat{e}$ and adding back the edges $e$ and $e'$.  By assumption $\hat{e}$ is contained in some $\hat{G}_j$ and we modify $\hat{G}_j$ by deleting $\hat{e}$ and adding the edge $e'$.  Finally, we add the edge $e$ to the graph $\hat{G}_i$ for $i = \phi(e)$.
We claim that after reversing the split of all $B$ edges, the graph $\hat{G}_i$ doesn't have $w$ as a cut vertex.
By definition, $i \in \hat{\phi}(\hat{A}\cup \hat{B})$. Now, if $i \in J$, then by the previous step (of reversing the split of $A$ edges) $e$ is not the only edge of $\hat{G}_i$ incident with $v$, so $w$ is not a cut vertex of $\hat{G}_i$. If $i \notin J$, then $i \in \hat{\phi}(\hat{B})$,
say $i=\hat{\phi}(f)$ for some $f \in \hat{A}$. Then, after reversing the split of $f$, we again get that $w$ is not a cut vertex of $\hat{G}_i$.

Let us pause the action here, after all splits involving edges of $A \cup B$ have been reversed.  At this point the graphs
$\hat{G}_1, \ldots, \hat{G}_k$ satisfy properties (C1) and (C3).  Also, property (C4) holds for every vertex in $T \setminus \{w,v\}$.  At present, the vertex $v$ has degree at most $2k$, but it may be that all edges incident with it appear in the graphs
$\hat{G}_1, \ldots, \hat{G}_k$.  To arrange for property (C4) to hold at the vertex $v$ we will arrange for some unused edges at $v$ in the next part of the process.

\begin{figure}[h]
  \centering
  \includegraphics[width=200pt]{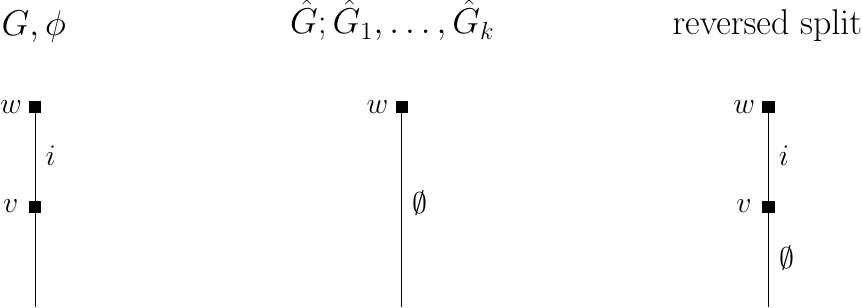}
  \caption{reversing the split of a $C$ edge.}
  \label{bbfig4}
\end{figure}

Next we consider reversing the splits of edges in $C$ (see Figure \ref{bbfig4}).
For every edge $e \in C$, note that $e$ was split off with some other edge $e'$ to form an edge $\hat{e} \in \hat{C}$.  By construction $\hat{e}$ will not be contained in any $\hat{G}_i$ and we now modify $\hat{G}$ by deleting $\hat{e}$, and adding back the edges $e$ and $e'$.  We add the edge $e$ to the graph $\hat{G}_i$ for $\phi(e) = i$.  We do not add the edge $e'$ to any $\hat{G}_i$.

\begin{figure}[h]
  \centering
  \includegraphics[width=200pt]{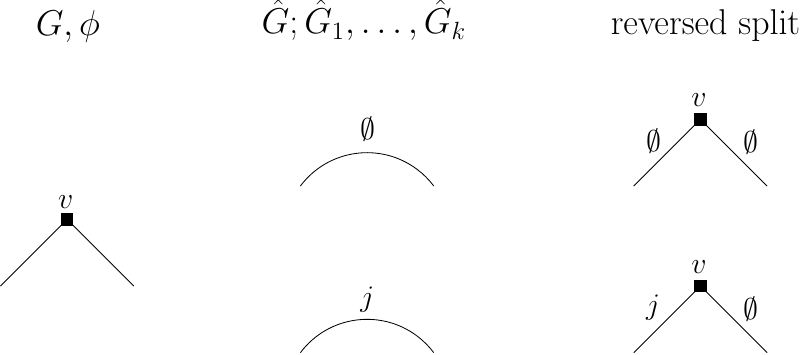}
  \caption{reversing the split of the remaining edges.}
  \label{bbfig5}
\end{figure}

Finally, we will reverse the splits of the remaining edges to return to the original graph (see Figure \ref{bbfig5}).
Every remaining edge $\hat{e}$ which is in our modified $\hat{G}$ but not $G$ was formed by splitting off a pair of edges $e, e'$ incident with $v$ with the property that neither $e$ nor $e'$ is incident with $w$.  We reverse this operation by deleting $\hat{e}$ and adding back the edges $e$ and $e'$ to the graph $\hat{G}$.  If $\hat{e}$ is not contained in any of the $\hat{G}_i$ graphs then we do not modify any of these graphs.  Next suppose that $\hat{e}$ is in contained in $\hat{G}_j$.  Temporarily let us modify $\hat{G}_j$ by deleting the edge $\hat{e}$ and adding both edges $e$ and $e'$.  After this operation, the graph $\hat{G}_j - w$ will contain a cycle using either $e$ or $e'$ (since $\hat{G}_j$ already spanned $v$).  Accordingly, we may delete one of $e$ or $e'$ from the graph $\hat{G}_j$ and still maintain properties (C1) and (C3).

When the process terminates, the graph $\hat{G}$ will be equal to our original graph $G$.  Furthermore, the vertex $v$ will be incident with at least $\frac{1}{2}( f(k) - 2k )\geq k$ edges which do not appear in any of the $\hat{G}_i$ graphs.  It follows that these graphs $\hat{G}_1, \ldots, \hat{G}_k$ satisfy properties (C1)-(C4) thus giving a solution to the original problem.  This contradicts the assumption that $G$ is a counterexample, completing the proof.  \hfill$\Box$

\bigskip

Now  we shall, as promised, reduce our original problem for Lemma \ref{steiner_lem} to a new one on a smaller graph.  To this end, we define the graph $G'$ to be the graph obtained from $G$ by deleting $w$ and every neighbour of $w$ contained in $V(G) \setminus T$ and then adding the edge set $D^{\circ}$ (see Figure \ref{figSpecialVertex}). The edges $D^{\circ}$ shall, in a way, represent the edges of $W_n^{\circ}$ that we removed. While we won't add any edges to $G'$ to represent the other edges we deleted, we keep track of those from $W_n^*$ by defining a set of edges (not appearing in $G'$)  $D^*_i$ for every $1 \le i \le k$ as follows.  If there exists $wv \in W^*_n$ with $\phi(wv) = i$ then we let $D^*_i = \{ v' v'' \}$, where $v' v'' \in D^*$ is the edge corresponding to $wv$ (so $N(v) = \{w,v',v''\}$)
and otherwise we set $D^*_i = \emptyset$.  So, the collection of all nonempty sets $D^*_i$ forms a partition of $D^*$.

We will use the previously defined values $p(v), q(v), u(v)$ to keep track of the number of deleted edges with nonzero and zero colours which were deleted in moving from $G$ to $G'$.  Note that our definitions imply that every $v \in T'$ satisfies
\begin{equation}
\label{gpdeg}
{\mathit deg}_{G'}(v) = f(k) - p(v) - q(v) - u(v).
\end{equation}
In the final solution every vertex in $T$ must be incident with at least $k$ edges which are not in any $G_i$ graph.  In order to keep track of the number of edges in $G'$ which we will be permitted to use at a vertex, we introduce the following capacity function.
\begin{equation}
\label{capdef}
{\mathit cap}(v) = f(k) - k - 2 - p(v) - q(v).
\end{equation}
The reason for including the $-2$ will be made apparent shortly.  For now let us check that this function is nonnegative.

\begin{observation}
\label{cappos}
${\mathit cap}(v) \ge 0$ for every $v \in T'$.
\end{observation}

\noindent{\it Proof:} By Observation \ref{obp2u} and Lemma \ref{buzzer_beater} we find
\begin{align*}
2 {\mathit cap}(v)
	&= 2f(k) - 2k - 4 - 2 p(v) - 2q(v)	\\
	&\ge 2f(k) - 2k - 4 - 2(k-1) - (f(k) + 1) \\
	&= f(k) -4k -3
	\ge 0. \quad\quad\Box
\end{align*}

Our next lemma will be used to reduce our problem on the original graph $G$ to a new one on the graph $G'$.

\begin{lemma}
\label{prime}
There do not exist edge-disjoint subgraphs $G_1',\ldots,G_k'$ of $G'$ satisfying all of the following conditions:
\begin{itemize}
\item[(D1)] $G_i'$ spans $T'$ for every $1 \le i \le k$.
\item[(D2)] $G_i' \cup D^*_i$ is connected for every $1 \le i \le k$.
\item[(D3)] Every $v \in T'$ satisfies $\sum_{i=1}^k {\mathit deg}_{G'_i}(v) \le cap(v)+2$.
\end{itemize}
\end{lemma}

\noindent{\it Proof:}
Suppose (for a contradiction) that there exist such graphs $G_1',\ldots,G_k'$.  Then we modify each $G_i'$ to obtain a subgraph $G_i$ of $G$ as follows.  For every edge $v'v'' \in D^{\circ} \cap E(G_i')$ let $wv$ be the corresponding edge of $W^{\circ}_n$ and add the vertex $v$ and the edges $vv'$ and $vv''$ to $G_i$.  Then, add all edges in $\phi'^{-1}(i)$ together with their endpoints (recall our definition of $\phi'$ from earlier).  The resulting graphs
$G_1,\ldots,G_k$ clearly satisfy (C1) and (C2). They satisfy (C3) because any edge in $D^{\circ}\cap E(G_i)$ is replaced by a 2-path of colour $i$ in $G_i$, and because we have (D2) and any edge in $D^*_i$ is replaced by a 2-path of colour $i$ in $G_i$ (by definition of $\phi'$). Condition (C4) holds because of (D3).
\hfill$\Box$

\bigskip

Next we will introduce some hypergraphs so as to take advantage of the matroidal properties of spanning hypertrees.  To do this, we define $H$ to be the hypergraph on vertex set $T$ obtained from $G$ by replacing each vertex $v \in V(G) \setminus T$ with an edge of size 3 containing $N(v)$.  Similarly, we let $H'$ denote the hypergraph obtained from $G'$ by replacing each $v \in V(G') \setminus T'$ with an edge of size 3 containing $N(v)$.  Note that by this definition and equation (\ref{gpdeg}), every $v \in T'$ satisfies
\begin{equation}
\label{hpdeg}
{\mathit deg}_{H'}(v) = f(k) - p(v) - q(v) - u(v).
\end{equation}
For simplicity of notation, let $E=E(H)$ and $E'=E(H')$.
It follows from the connectivity properties of $G$ that the hypergraph $H$ will be $f(k)$-edge-connected, so $H'$ will not have too many small edge-cuts.

\begin{lemma}
\label{edge_count}
If ${\mathcal P}'$ is a partition of $T'$ then
\[ 3 \lambda^{out}_{\mathcal P'}(E') \ge f(k) (|{\mathcal P'}| - 1)  -  \lambda^{out}_{\mathcal P'}(D^*) - f(k) + 2k \]
and if $|{\mathcal P}'| = 2$ this may be improved by $k - \tfrac{1}{2} \lambda^{out}_{\mathcal P'} (D^*)$ to
\[  3 \lambda^{out}_{\mathcal P'}(E') \ge 3k - \tfrac{3}{2} \lambda^{out}_{\mathcal P'}(D^* ). \]
\end{lemma}

\noindent{\it Proof:}
Extend ${\mathcal P}'$ to a partition ${\mathcal P}$ of $T$ by adding $\{w\}$ as a block.
Then ${\mathcal P}$ is a partition of $V(H)$ and since $H$ is $f(k)$-edge-connected, we have that every $B \in {\mathcal P}$ satisfies $d_H(B)=|\delta_H(B)| \ge f(k)$.
Let $W^*$ be the set of hyperedges in $H$ associated with an edge in $W^*_n \cup W^*_t$ and similarly let $W^{\circ}$ be the set of hyperedges in $H$ associated with an edge in $W^{\circ}_n \cup W^{\circ}_t$.  For $i=2,3$ let $W^*_i$
 $(W^{\circ}_i)$ denote the set of hyperedges in
$W^*$ $(W^{\circ})$ that have ends in exactly $i$ different blocks of ${\mathcal P}$.
Similarly, for $i=2,3$, we define $E_i$ to be the set of all edges not in $W^* \cup W^{\circ}$ which have ends in exactly $i$ different blocks of ${\mathcal P}$. Now every edge in $W_3^{\circ}$ corresponds to an edge in $D^{\circ}$ which has ends in distinct blocks of ${\mathcal P}$ and similarly every edge in $W_3^*$ corresponds to an edge in $\Lambda_{\mathcal P'}^{out}(D^*)$ (though the edges in $D^*$ are not in $H'$).
We now have
\begin{align*}
	3 \lambda^{out}_{\mathcal P'}(E')	
	&= 3 \lambda^{out}_{\mathcal P}(E)  + 3 |W^{\circ}_3 | - 3 f(k)\\
	&= 	\left(	 \sum_{B \in {\mathcal P}} d_{H}(B) + |E_2| + |W_2^*|  +  |W_2^{\circ}|\right)	+ 3 |W^{\circ}_3 | - 3f(k)\\
	&\ge	\left(	f(k) |{\mathcal P}|  + |W_2^*| + |W_2^{\circ}| \right) + |W^{\circ}_3 | - 3 f(k).
	\end{align*}	
Note in the second line that $\sum_{B\in \mathcal{P}}d_H(B)$ counts hyperedges with ends in 3 blocks of $\mathcal{P}$ three times, and those with ends in 2 blocks of $\mathcal{P}$ only twice.
Now, using $|W_2^*|+|W_3^*|=k$ and  $|W_2^{\circ}|+|W_3^{\circ}|=k$ (by (P1)), this gives us
	\begin{align*}			
    &= 	f(k) \left(|{\mathcal P'}|+1\right)  + \left( k- |W_3^*| \right) + k - 3 f(k)	\\
	&=		f(k) (|{\mathcal P}'| - 1)  - |W^*_3| - f(k) + 2k 	\\
	&=		f(k) (|{\mathcal P}'| - 1)  - \lambda^{out}_{\mathcal P'} (D^* ) - f(k) + 2k,
\end{align*}
as desired.  If $|{\mathcal P}'| = 2$ then we set $F$ to be the set of edges in $H$ which do not contain $w$
but do intersect both blocks of ${\mathcal P}'$ (note that $F\subseteq E_2$, but $E_2$ could contain edges coloured by $\phi$ incident to $w$.) We now repeat the same calculation.
\begin{align*}
	3 \lambda^{out}_{\mathcal P'}(E')	
	&= 3 \lambda^{out}_{\mathcal P}(E)  + 3 |W^{\circ}_3 | - 3 f(k)	\\
	&= 		\left( \sum_{B \in {\mathcal P}} d_{H}(B) + |E_2| + |W_2^*|  +  |W_2^{\circ}|\right)	+ 3 |W^{\circ}_3 | - 3f(k)	\\
	&\ge		3 f(k) + |F| + |W_2^*| + |W_2^{\circ}|  + 3|W^{\circ}_3 | - 3 f(k)	\\
	&\ge		\Big( |F| + |W^{\circ}_3| \Big) + |W_2^*| + \Big( |W^{\circ}_2| + |W^{\circ}_3| \Big)	 \\
&=\Big( |F| + |W^{\circ}_3| \Big) +\Big(k-|W^*_3|\Big) + k	 \\
	&=		\lambda^{out}_{\mathcal P'} (E' ) + 2k  - \lambda^{out}_{\mathcal P'} (D^* ).
\end{align*}
 Rearranging the final inequality we get that $3 \lambda^{out}_{\mathcal P'}(E') \ge 3k - \tfrac{3}{2} \lambda^{out}_{\mathcal P'} (D^* )$ as desired.
\hfill$\Box$

\bigskip

With this, we are finally ready to prove Lemma \ref{steiner_lem}.

\bigskip

\noindent{\it Proof of Lemma \ref{steiner_lem}:} We will establish this lemma by contradicting Lemma \ref{prime}, thus contradicting the assumption that $G$ is a counterexample.  As suggested earlier, our proof will make use of hypergraphic matroids. We have already defined the hypergraph $H'$ and $E' = E(H')$. We further define $H''$ to be the hypergraph obtained from $H'$ by adding the edges $D^*$, and set $E'' = E(H'')$. Finally, we let $M$ be the hypergraphic matroid on $H''$, and let $r$ be the rank function of $M$.

 We claim that there cannot exist a vector $x \in {\mathbb R}_+^{E'}$ satisfying all of the following properties:
 (note that in the following equation, the rank function is based on the hypergraphic matroid $H''$ though many equations involve the hypergraph $H'$)
\begin{itemize}
\item $x(e) \le 1$ for every $e \in E'$
\item $x( \Lambda^{in}_{\mathcal P}( E') )  \le k r( {\mathcal P} ) - \lambda^{in}_{\mathcal P}(D^* )$ for every partition ${\mathcal P}$ of $T'$
\item $-x(E') \le - k r(M) + |D^*|$
\item $x( \delta_{H'}(v) ) \le {\mathit cap}(v)$ for every $v \in T'$.
\end{itemize}
To establish the claim, suppose, for a contradiction, that such an $x$ does exist. Define the vector $\bar{x} \in {\mathbb R}^{E''}$
by the rule that $\bar{x}(e) = x(e)$ for every $e \in E'$ and $\bar{x}(e) = 1$ for every $e \in D^*$.
For every partition $\mathcal{P}$ of $T'$,
$$\bar{x}(\Lambda^{in}_{\mathcal{P}}(E''))= x( \Lambda^{in}_{\mathcal P}( E') )+ \lambda^{in}_{\mathcal P}(D^* )\leq k\, r(\mathcal{P}).$$
When $\mathcal{P}$ is the partition with only one block, the second condition becomes
$$x(E') \leq k(|T'|-1)-|D^*|=k r(M)-|D^*|,$$
which gives
$$\bar{x}(E'')=x(E')+|D^*|=k\,r(M)=k (|T'|-1).$$

Therefore Lemma \ref{partition_rank} tells us that $\bar{x}$ is a fractional basis of $kM$ and $kM$ has rank $k\,r(M)$.

Since $\bar{x}(e) = 1$ for every $e \in D^*$, the vector $x$ is a fractional basis of the matroid $L = (kM)/D^*$ on $E'$.  We apply Theorem \ref{kls} to the matroid $L$, the vector $x$, and the family of subsets of $E'$ given by $\{ \delta_{H'} (v) \mid v \in T' \}$.  Since every edge is contained in at most three of these sets, the theorem gives us a basis $B$ of $L$ for which
\begin{equation}\label{C3}
|B \cap \delta_{H'}(v)| \le cap(v)+2
\end{equation}
holds for every $v \in T'$.  Since $B \cup D^*$ is a basis of $kM$ and $|D^*| \le k$, Lemma \ref{adjust_union} says that we may choose a partition of $B \cup D^*$ into bases $\{B_1, \ldots B_k \}$ of $M$ so that $D^*_i \subseteq B_i$ for every $1 \le i \le k$.  Now, for each $1 \le i \le k$ we let $G_i'$ denote the subgraph of $G'$ which is constructed as follows.  We begin with $G_i'$ containing the vertex set $T'$ and no edges.  Then we add to $G_i'$ all edges of size 2 in $B_i \setminus D^*_i$, and then for every edge of size 3 in
$B_i \setminus D_i^*$ containing $\{v,v',v''\}$ this edge corresponds to a vertex $x \in V(G') \setminus T'$ in $G'$ and we add to $G'_i$ the vertex $x$ and the edges $xv, xv', xv''$.
Since $B_i$ is a basis of $L$, $G_i'$ is connected and spans $T'$.
Therefore each $G_i'$ satisfies (D1), (D2), and (D3). Thus, these graphs contradict  Lemma \ref{prime}, and our claim has been established.

Since we know that $x$ cannot exist, Farkas' Lemma tells us that there must exist families of nonnegative parameters $\{y_e\}_{e \in E'}$, $y_G$, $\{ y_v \}_{v \in T'}$, and $\{ y_{\mathcal P} \}$ for every partition ${\mathcal P}$ of $T'$, which satisfy the following properties
\begin{equation}
\label{eqa1}
y_G \Big( k (|T'|-1) - |D^*| \Big) > \sum_{e \in E'} y_e + \sum_{ {\mathcal P} } y_{\mathcal P} \Big( k \, r({\mathcal P}) - \lambda^{in}_{\mathcal P}(D^* ) \Big) + \sum_{v \in T'} y_v \, {\mathit cap}(v)
\end{equation}
\begin{equation}
\label{eqb1}
y_e + \sum_{ {\mathcal P} : e \in \Lambda^{in}_{\mathcal P}(E') } y_{\mathcal P} + \sum_{v \sim e} y_v \ge y_G \mbox{\quad\quad for every $e \in E'$}
\end{equation}

We claim that we may assume $y_G=1$. To this end, it suffices to show that $y_G>0$, as in this case we may scale $y$ by a positive constant.   To prove this, we will consider equation (\ref{eqa1}).  First note that Observation \ref{cappos} shows that ${\mathit cap}(v) \ge 0$ for every $v \in T'$.  Since $\lambda^{in}_{\mathcal P}(D^* ) \le |D^*| \le k$ for every partition, we also know that the coefficient of $y_{\mathcal{P}}$ in (\ref{eqa1}) is  nonnegative.  Finally, $|D^*|\leq k$ (by definition) and $|T'|\geq 2$ (by (P5)), so we must indeed have $y_G>0$, and hence can assume $y_G=1$.

Next we claim that we can choose $y$ so that, if ${\mathcal P}_1$ and ${\mathcal P}_2$ are partitions with $y_{ {\mathcal P}_1 }, y_{ {\mathcal P}_2 } > 0$, then it must be that they are nested, i.e., either ${\mathcal P}_1 \le {\mathcal P}_2$ or ${\mathcal P}_2 \le {\mathcal P}_1$. To do so, over all vectors $y$ which satisfy equations (\ref{eqa1}), (\ref{eqb1}), and $y_G=1$ we choose $y$ to be extreme in the following sense
\begin{enumerate}
\item $\sum_{ {\mathcal P} } y_{\mathcal P} r({\mathcal P}) $ is minimum.
\item $\sum_{ {\mathcal P} } y_{\mathcal P} (r( {\mathcal P} ))^2$ is maximum (subject to T1).
\end{enumerate}

Now suppose that ${\mathcal P}_1$ and ${\mathcal P}_2$ are partitions with $y_{ {\mathcal P}_1 }, y_{ {\mathcal P}_2 } \ge \epsilon > 0$,
which  are not nested. Then consider modifying our vector $y$ by decreasing $y_{ {\mathcal P}_1 }$ and $y_{ {\mathcal P}_2 }$ by $\epsilon$ and then increasing $y_{ {\mathcal P}_1 \vee {\mathcal P}_2 }$ and
$y_{ {\mathcal P}_1 \wedge {\mathcal P}_2 }$ by $\epsilon$.
 It follows immediately that this new $y$ still satisfies equation (\ref{eqb1}). Similarly, note that
 $$\lambda^{in}_{{\mathcal P}_1 \vee {\mathcal P}_2} (D^*) + \lambda^{in}_{{\mathcal P}_1 \wedge {\mathcal P}_2}(D^* ) - \lambda^{in}_{ {\mathcal P}_1}(D^* ) - \lambda^{in}_{ {\mathcal P}_2}(D^* )=0,$$
 so the effect of modifying $y$ increases the right-hand side of (\ref{eqa1}) by
 $$\epsilon\, k\Big( r( {\mathcal P}_1 \vee {\mathcal P}_2 ) + r( {\mathcal P}_1 \wedge {\mathcal P}_2 ) - r( {\mathcal P}_1 ) - r( {\mathcal P}_2 ) \Big),$$
and by submodularity of the rank function $r$, this is nonpositive. This operation will strictly improve $y$ with respect to our first optimization criteria unless $r( {\mathcal P}_1 \vee {\mathcal P}_2 ) + r( {\mathcal P}_1 \wedge {\mathcal P}_2 ) = r( {\mathcal P}_1 ) + r( {\mathcal P}_2 )$.  However, in this case the second induction criteria is affected by an increase of
\[ \epsilon \Big(  \big( r( {\mathcal P}_1 \vee {\mathcal P}_2 ) \big)^2 +  \big( r( {\mathcal P}_1 \wedge {\mathcal P}_2 ) \big)^2 - \big( r( {\mathcal P}_1 ) \big)^2 - \big( r( {\mathcal P}_2 ) \big)^2 \Big) > 0. \]
Therefore, we conclude that we may choose $y$ so that the only partitions ${\mathcal P}$ with $y_{ \mathcal P} > 0$ are nested.

Next we claim we may assume that for the trivial partition $\{T'\}$, $y_{ \{T'\} } = 0$.  To see this, suppose that $y_{ \{T'\} } = \delta > 0$ and consider modifying $y$ by decreasing both $y_G$ and $y_{ \{T'\} }$ by $\delta$. Then the modified $y$ is still nonnegative and still satisfies
(\ref{eqb1}). Moreover, it also still satisfies (\ref{eqa1}) because
  $k\,r(\{T'\}) - \lambda^{in}_{\{T'\}}(D^* )=k (|T'|-1) - |D^*|$.
Hence we can choose $y$ so that $y_{T'}=0$. After this last step, we may again rescale $y_G$ and all other entries of $y$ to get $y_G=1$ again, and the property on nested partitions still holds.

Since the trivial partition $\{T'\}$ has zero weight, we may assume that there are partitions ${\mathcal P}_0 \le {\mathcal P}_1 \le \ldots \le {\mathcal P}_{\ell}$ so that every partition ${\mathcal Q} \not\in \{ {\mathcal P}_0, \ldots, {\mathcal P}_{\ell} \}$ has $y_{\mathcal Q} = 0$ and so that
${\mathcal P}_0$ is the discrete partition and ${\mathcal P}_{\ell}$ is a partition of size 2.  For convenience of notation, we now set $y_i = y_{ {\mathcal P}_i }$, $\Lambda^{out}_i(F) = \Lambda^{out}_{ {\mathcal P}_i }(F)$, $\Lambda^{in}_i(F) = \Lambda^{in}_{ {\mathcal P}_i }(F)$ and assign $\lambda_i^{out}$ and $\lambda^{in}_i$ similarly for every $0 \le i \le \ell$ and every subset of edges $F$.  With this, we now restate equation (\ref{eqb1}) as follows:
\begin{equation}
\label{eqb2}
y_e + \sum_{ i \in \{0,\ldots, \ell\} : e \in \Lambda_i^{in}(E') } y_i + \sum_{v \sim e} y_v \ge 1 \mbox{ \hspace*{.2in} for every $e \in E'$}.
\end{equation}

Next we show that we can further modify $y$, preserving all the properties discussed so far, to get the additional property that $\sum_{i=1}^{\ell} y_i=1$.
If $\sum_{i=1}^{\ell} y_i > 1$ then let $j$ be the minimum over all $\{1,\ldots,\ell\}$ for which $y_j > 0$.  Then, in Equation (\ref{eqb2}) every edge in $\Lambda^{in}_j(E' )$ receives a total weight of at least $\sum_{i=1}^{\ell} y_i > 1$ from the partitions so none of these edges are tight with respect to this constraint, and then we may decrease $y_j$ by $\min\{ y_j, \sum_{i=1}^{\ell} y_i - 1 \}$. Possibly repeating this process (if $\min\{ y_j, \sum_{i=1}^{\ell} y_i - 1 \}=y_j$),  we can get
$\sum_{i=1}^{\ell} y_i \le 1$.  Now, the term $y_0$ associated with the discrete partition has no effect on any of our equations.  Therefore, we may choose $y_0$ so that $\sum_{i=0}^{\ell} y_i = 1$.  This normalization is handy and we use it next to rewrite the following terms:
\begin{align*}
&|D^*| - \sum_{i=0}^{\ell} y_i \, \lambda^{in}_i(D^*) = \sum_{i=0}^{\ell} y_i \big( |D^*| - \lambda^{in}_i(D^* ) \big) = \sum_{i=0}^{\ell} y_i \, \lambda_i^{out}(D^* ),\\
\textrm{and}\\
&k( |T'| - 1 ) - \sum_{i=0}^{\ell} y_i k r( {\mathcal P}_i) = \sum_{i=0}^{\ell} y_i k \big( |T'| - 1 - r ( {\mathcal P_i} ) \big) = \sum_{i=0}^{\ell} y_i k ( |{\mathcal P}_i| - 1 ).
\end{align*}
Taking advantage of the above identities we may now rewrite equation (\ref{eqa1}) as follows
\begin{equation}
\label{eqa2}
\sum_{i=0}^{\ell}  y_i \Big( k( |{\mathcal P}_i| - 1) - \lambda^{out}_i(D^*) \Big)  > \sum_{e \in E'} y_e + \sum_{v \in T'} y_v \, {\mathit cap}(v)
\end{equation}

In the remainder of our proof we will use (\ref{eqa2}) and (\ref{eqb2}), along with some previous results, to achieve our desired contradiction. In the following calculation, the first inequality is obtained by summing (\ref{eqb2}) over all edges, the second inequality uses Lemma \ref{edge_count} (including the special case that $|{\mathcal P}_{\ell}| = 2$), and the final inequality uses equation (\ref{eqa2}):
\begin{align*}
 3 \sum_{e \in E'} y_e &+ 3\sum_{v \in T'} y_v {\mathit deg_{H'}}(v)
	\ge 3|E'| -  3\sum_{j=0}^{\ell} y_j  \, \lambda^{in}_j (E') 	\\
	&= 3\sum_{j=0}^{\ell}  y_{j} \, \lambda^{out}_j(E') 	\\
	&\ge \sum_{j=0}^{\ell}   y_j \Big( f(k)( |{\mathcal P}_j| - 1) - \lambda^{out}_j(D^*) 		
			- f(k) + 2k \Big)   + y_{\ell} (k - \tfrac{1}{2} \lambda^{out}_{\ell}(D^* )) \\
	&\ge \tfrac{f(k)}{k}  \sum_{j=0}^{\ell}   y_j \Big( k( |{\mathcal P}_j| - 1) - \lambda^{out}_j(D^*) \Big) 		
			-  f(k) + 2k   + y_{\ell} \big( k + \tfrac{f(k)-3k}{2k} \lambda^{out}_{\ell}(D^* ) \big)\\
	&> \tfrac{f(k)}{k} \sum_{e \in E'} y_e + \tfrac{f(k)}{k}  \sum_{v \in T'} y_v  {\mathit cap}(v)
		- f(k) + 2k + y_{\ell} \big( k + \tfrac{f(k)-3k}{2k} \lambda^{out}_{\ell}(D^* ) \big).
\end{align*}
Rearranging brings us to the following inequality
\begin{equation}\label{eqfk2} f(k) - 2k > \tfrac{1}{k} \sum_{v \in T'} y_v \Big( f(k) {\mathit cap}(v) - 3k {\mathit deg}_{H'}(v) \Big) + \tfrac{f(k) - 3k}{k}  \sum_{e \in E'} y_e
	+ y_{\ell} \big( k + \tfrac{f(k)-3k}{2k} \lambda^{out}_{\ell}(D^* ) \big) .\end{equation}

In the next calculation we consider an arbitrary vertex $v \in T'$.  We call upon the degree formula ${\mathit deg}_{H'}(v) = f(k) - p(v) -q(v) - u(v)$ (equation (\ref{hpdeg})) and the definition of the capacity function, ${\mathit cap}(v) = f(k) - k - 2 - p(v) - q(v)$ (equation (\ref{capdef})); the inequalities $q(v) \le u(v) + 1$ (Observation \ref{obp2u}) and $p(v) \le k-1$ (Lemma \ref{buzzer_beater}); as well as the basic bounds $k \ge 3$ and  $5k+3 \le f(k) \le 5k+4$.
\begin{align*}
&f(k) {\mathit cap}(v) - 3k {\mathit deg}_{H'}(v)\\
	&\ge f(k) (4k+1-p(v) - q(v))  - 3k (f(k) - p(v) - q(v) - u(v))	\\
	&= 	(k+1) f(k) - (f(k) - 3k) (p(v)+q(v))   + 3k u(v)	\\
	&\ge	(k+1)f(k) - (f(k) - 3k)(k-1) - (f(k) - 3k)	\\
	&=	3k^2 + f(k) \\
	&\ge k(f(k) - 2k).
\end{align*}
Combining this with (\ref{eqfk2}) gives us the following:
\begin{equation}\label{fk2k} (f(k) - 2k) \Big( 1 - \sum_{v \in T'} y_v \Big) > \tfrac{f(k) - 3k}{k} \sum_{e \in E'} y_e  + y_{\ell} \big( k + \tfrac{f(k)-3k}{2k}  \lambda^{out}_{\ell}(D^* ) \big). \end{equation}

Setting $\alpha = 1 - \sum_{v \in T'} y_v$ it follows from (\ref{fk2k}) that $\alpha > 0$.  Furthermore, every edge $e \in \Lambda^{out}_{\ell}(E')$ must satisfy $y_e \ge \alpha$ by (\ref{eqb2}), since this edge is outer with respect to every partition ${\mathcal P}_i$ and $\sum_{v \sim e} y_v \le \sum_{v \in T'} y_v = 1- \alpha$.  We are now able to derive a contradiction based on the value of $y_{\ell}$. First suppose that $y_{\ell} \ge \alpha$.  In this case (\ref{fk2k}) gives us the following contradiction (here we have applied Lemma \ref{edge_count} for the last inequality).
\begin{align*}
(f(k) - 2k) \alpha
	&> \tfrac{f(k) - 3k}{k} \sum_{e \in E'} y_e + y_{\ell} \big( k + \tfrac{f(k)-3k}{2k}  \lambda^{out}_{\ell}(D^* ) \big) 	\\
	&\ge \tfrac{f(k)-3k}{k} \alpha \,\lambda^{out}_{\ell}(E' )  + \alpha  \big( k + \tfrac{f(k)-3k}{2k} \lambda^{out}_{\ell}(D^* ) \big)	\\
	&= \tfrac{f(k)-3k}{k} \alpha \Big( \lambda^{out}_{\ell}(E') + \tfrac{1}{2}  \lambda^{out}_{\ell}(D^* ) \Big) + \alpha k	\\
	&\ge (f(k) - 2k) \alpha.
\end{align*}
It follows that $\beta = \alpha - y_{\ell} \ge 0$.  Now, every edge in $\Lambda_{\ell-1}(E') \setminus \Lambda_{\ell}(E' )$ must satisfy
$y_e \ge \beta$  by (\ref{eqb2}), since this edge receives a total contribution of at most $1- \alpha$ from the vertices and $y_{\ell}$ from the partitions.  Hence (\ref{fk2k}) gives us
\[
(f(k) - 2k) \alpha >  \tfrac{f(k) - 3k}{k} (\alpha - \beta) \lambda^{out}_{\ell}(E') + \tfrac{f(k)-3k}{k} \beta \, \lambda^{out}_{\ell-1}(E') + (\alpha - \beta) \big( k + \tfrac{f(k)-3k}{2k} \lambda^{out}_{\ell}(D^* ) \big).
\]
Splitting the left hand side into $(f(k) - 2k) (\alpha - \beta) + (f(k) - 2k) \beta$ and then rearranging gives:
\begin{equation}\label{alphabeta} 0 > (\alpha - \beta) \Big( \tfrac{f(k) - 3k}{k} \big(\lambda^{out}_{\ell}(E' ) + \tfrac{1}{2} \lambda^{out}_{\ell}(D^*) \big) + k  - (f(k) - 2k) \Big)
	+ \beta \Big( \tfrac{f(k)-3k}{k} \lambda^{out}_{\ell-1}(E')  - (f(k) - 2k) \Big). \end{equation}
As in the previous case the term above with coefficient $\alpha - \beta$ must be $\ge 0$ by Lemma \ref{edge_count}.  Using the same lemma, we find that
$$\lambda^{out}_{\ell-1}(E') \ge \tfrac{1}{3} ( f(k) + 2k - \lambda^{out}_{\ell -1}(D^*) ) \ge \tfrac{1}{3}(f(k) + k) \ge 2k$$
and since the coefficient of $\beta$ in (\ref{alphabeta}) must be negative we have
\[ 0 > \big( \tfrac{f(k)-3k}{k}(2k) - ( f(k) - 2k ) \big) = f(k) - 4k \ge 0,  \]
and this final contradiction completes the proof.
\hfill$\Box$

\section{Packing Connectors}

The proof of Theorem \ref{connector_thm} is very similar to that of Theorem \ref{steiner_thm} so we shall not give this proof in full detail, but rather we step through the argument referring heavily to the previous section. Our main lemma will now involve the function $g(k) =6k+6$.

\begin{lemma}
\label{connector_lem}
Let $G$ be a graph and let $T \subseteq V(G)$ satisfy $|T| \ge 2$ and have the property that every edge cut separating $T$ has size $\ge g(k)$.  Let $w \in T$ be a distinguished vertex with ${\mathit deg}(w) = g(k)$ and let $\mu : \delta(w) \rightarrow \{0,1,\ldots,k\}$
be onto and satisfy $| \mu^{-1}(0) | \ge 2k$.  Then there exist $k$ edge-disjoint subgraphs $G_1,\ldots,G_k$ which satisfy
the following conditions:
\begin{itemize}
\item[(E1)] $G_i$ is a $T$-connector for $1 \le i \le k$.
\item[(E2)] $\delta(w) \cap E(G_i) = \mu^{-1}(i)$ for $1 \le i \le k$.
\item[(E3)] $w$ is not a cut vertex of $G_i$ for $1 \le i \le k$.
\item[(E4)] Every $v \in T$ satisfies $\sum_{i=1}^k {\mathit deg}_{G_i}(v) \le deg_G(v) - 2k$.
\end{itemize}
\end{lemma}

The proof of Theorem~\ref{connector_thm} follows from Lemma~\ref{connector_lem} just as
Theorem~\ref{steiner_thm} follows from Lemma~\ref{steiner_lem}, and we therefore omit these details.
To prove Lemma \ref{connector_lem} we assume that $G$ together with $T$, $w$ and $\mu$ is a counterexample which is extreme in the following sense.
\begin{enumerate}
\item[(Y1)] $\sum_{v \in T} {\mathit deg}(v)$ is minimum.
\item[(Y2)] $|V(G)| + |E(G)|$ is minimum (subject to Y1).
\item[(Y3)] $|\mu^{-1}(0)|$ is minimum (subject to Y1 and Y2).
\end{enumerate}

By arguments analogous to those in Lemma \ref{minimal_counter} we have the following result.
\begin{lemma}
The graph $G$ satisfies the following properties:
\begin{itemize}
\item[(Q1)] $|\mu^{-1}(0)| = 2k$.
\item[(Q2)] Every $v \in V(G) \setminus T$ has degree 3 and has $|N(v)| = 3$.
\item[(Q3)] Every $v \in T$ has degree $g(k)$.
\item[(Q4)] $V(G) \setminus T$ is an independent set.
\item[(Q5)] $|T|\geq 3$.
\end{itemize}
\end{lemma}

Condition (Q2) implies that our connectors will be connected subgraphs of $G$ which span $T$ and have the property that every vertex in $V(G)\setminus T$ has degree $0$ or $2$.
The fact that we are not allowed to use all three edges incident with a vertex in $V(G)\setminus T$ is the main difference from dealing with Steiner trees. Because of this, we will not distinguish edges $W_n^*$ (or $W_t^*$) as in Section 4. We still define $W^{\circ}_t$ to be the set of edges in $\mu^{-1}(0)$ with both ends in $T$ and let $W^{\circ}_n = \mu^{-1}(0) \setminus W^{\circ}_t$. However, now we set $W^{\bullet}_t$ to be the set of all edges in $\mu^{-1}( \{1,\ldots,k\} )$ with both ends in $T$ and set $W^{\bullet}_n = \mu^{-1}( \{1,\ldots,k\}) \setminus W^{\bullet}_t$.
As before, for each choice of $\star \in \{ \circ, \bullet \}$, let $X^{\star}$ be the set of endpoints of edges of the form $W^{\star}_n$ other than $w$.
We construct two sets of auxiliary edges $D^{\circ}$ and $D^{\bullet}$.  For each $\star \in \{ \circ, \bullet \}$ we construct $D^{\star}$ from $\emptyset$ as follows: for each $x \in X^{\star}$ let $\{w,v,v'\}$ be the neighbours of $x$ and we add a new edge with ends $v$ and $v'$ to $D^{\star}$.  As before this construction yields a natural correspondence between $W^{\star}_n$ and $D^{\star}$.

Next we define an extension $\mu'$ of our precolouring.
We define $\mu'(e)=\mu(e)$ for all $e\in \delta(w)$ and also define $\mu'$ on all edges  of $\delta(x)$ where $x\in X^{\bullet}$.  We define $G^{\bullet}$ to be the graph $(V(G), D^{\bullet})$ and then choose an orientation of $G^{\bullet}$ so that every vertex $v$ satisfies
  \begin{equation}
 \label{orientation2}
 {\mathit deg}^+_{G^{\bullet}}(v) \le {\mathit deg}^-_{G^{\bullet}}(v) + 1.
 \end{equation}
Now for each $x\in X^{\bullet}$ with $N(x)=\{w,v,v'\}$, we consider the edge $vv' \in D^{\bullet}$.  If this edge is oriented from $v$ to $v'$ then we define
\begin{align*}
	\mu'(xv) &= \mu(wx)	\\
	\mu'(xv') &= 0
\end{align*}

This choice for the extended precolouring  has the property that every vertex in $X^{\bullet} \cup X^\circ$ has either only incident edges with a zero colour or exactly two incident edges precoloured with the same colour. Therefore the precolouring will not force us to use three edges incident to the same vertex in $V(G)\setminus T$ in a connector.
As for the Steiner tree case, we reduce the problem to one on a smaller graph obtained by deleting all edges precoloured by $\mu'$.
To keep track of the deleted edges we define the following functions for every $v \in T' = T \setminus \{w\}$
\begin{align*}
	p(v)	&= | \delta(v)\cap W^\bullet_t|\\
	q(v)	&= | \delta(v) \cap\mu'^{-1}(\{1,\ldots,k\})|-p(v)	\\
	u(v)	&= |\delta(v) \cap \mu'^{-1}(\{0\}) |.\\
\end{align*}
So $p(v) + q(v)$ is the number of edges at $v$ that are precoloured with a nonzero value by $\mu'$.
These edges have two types: their endpoint other than $v$ might be $w$, or might be in $X^{\bullet}$.  Here $p(v)$ counts the edges of the first type, and $q(v)$ counts the edges of the second type. The quantity $u(v)$ counts the number of edges incident with $v$ which are assigned the colour $0$ by $\mu'$ (the other end of such an edge could be either $w$ or a vertex in $X^{\bullet}$).

The following two results provide useful inequalities involving $p(v)$, $q(v)$ and $u(v)$ and are the analogs of Observation~\ref{obp2u} and Lemma~\ref{buzzer_beater} for connectors. The proof of Observarion~\ref{obp2uconn} is almost identical to the proof of Observation~\ref{obp2u} and therefore omitted.
The proof of Lemma~\ref{buzzer_beaterconn} is also omitted, as it is just a simplified version of the proof of Lemma~\ref{buzzer_beater} where the set $A$ is empty, and where in the end we have at least $\tfrac{1}{2}(g(k)-2k)\geq 2k$ edges which do not appear in any of the $\hat{G}_i$ graphs.
The only thing to notice for the proof of Lemma~\ref{buzzer_beaterconn} is that we may create final graphs $\hat{G}_1,\ldots,\hat{G}_k$ with vertices in $V(G)\setminus T$ of degree one. This can be easily corrected in the following way: if $\hat{G}_i$ has a vertex $u \in V(G)\setminus T$ of degree one, then we remove $u$ and the edge of $\hat{G}_i$ incident with $u$ to obtain a connector (by (E3) this edge was not precoloured).

\begin{observation}\label{obp2uconn}
Every $v\in T'$ satisfies
\begin{itemize}
\item $q(v) \le u(v) + 1$, and
\item $2 q(v) \le g(k) + 1$.
\end{itemize}
\end{observation}

\begin{lemma}
\label{buzzer_beaterconn}
Every $v \in T'$ satisfies $p(v) \le k-1$.
\end{lemma}

We need this extra result for dealing with connectors.

\begin{lemma}\label{p1plusp2}
Every $v\in T'$ satisfies $ p(v)+q(v)\leq 3k+3$.
\end{lemma}

\noindent{\it Proof:}
Let $v\in T'$, and consider the cut in $G$ induced by $w,v$, and all common neighbours of $w$ and $v$ that are in $X^{\bullet}$.
There are $q(v)+u(v)-|\delta_G(v)\cap W^\circ_t|$ of these common neighbours.
Since $|T|\geq 3$, this cut separates $T$ and hence must have size at least $g(k)$.
There are $p(v)+|\delta_G(v)\cap W^\circ_t|$ edges between $w$ and $v$, and both these vertices have degree $g(k)$, so the induced cut has size
\begin{align*}
& 2g(k) -2(p(v)+|\delta_G(v)\cap W^\circ_t|)-(q(v)+u(v)-|\delta_G(v)\cap W^\circ_t|)\\
&\le 2g(k) -2p(v)-q(v)-u(v)\\
&\le 2g(k) -2p(v)-2q(v)+1,
\end{align*}
where we used the first inequality in Observation~\ref{obp2uconn}.
Since the cut has size at least $g(k)$, we get
\begin{eqnarray*}
g(k)&\leq&2g(k) -2p(v)-2q(v)+1\\
\end{eqnarray*}
which implies that $p(v)+q(v) \leq \tfrac{g(k)+1}{2}$.
Now the result follows from the facts that $g(k)=6k+6$ and that the quantity $p(v)+q(v)$ is integral.\hfill$\Box$

We shall again reduce our original problem for Lemma \ref{connector_lem} to a new one on a smaller graph.  To this end, we define the graph $G'$ to be the graph obtained from $G$ by deleting $w$ and every neighbour of $w$ contained in $V(G) \setminus T$ and then adding the edge set $D^{\circ}$.
Our definitions imply that every $v \in T'$ satisfies
\begin{equation}
\label{gpdegconn}
{\mathit deg}_{G'}(v) = g(k) - p(v) - q(v) - u(v).
\end{equation}

In terms of the formulation of our new problem, a key quantity to keep track of, for each $v\in T'$, will be
\begin{equation}
\label{capdefconn}
{\mathit cap}(v) = g(k) - 2k - 2 - p(v) - q(v).
\end{equation}

One can show that this quantity is always nonnegative using Observation \ref{obp2uconn} and Lemma \ref{buzzer_beaterconn}
(this is similar to the proof of Observation~\ref{cappos}).

\begin{observation}
\label{capposconn}
${\mathit cap}(v) \ge 0$ for every $v \in T'$.
\end{observation}

The next result is the analog of Lemma~\ref{prime} for connectors. We omit the proof, as it is analogous to that of Lemma~\ref{prime} (with the omission of edges in $D^*$).

\begin{lemma}
\label{primeconn}
There do not exist edge-disjoint subgraphs $G_1',\ldots,G_k'$ of $G'$ satisfying all of the following conditions:
\begin{itemize}
\item[(F1)] $G_i'$ spans $T'$ for every $1 \le i \le k$.
\item[(F2)] $G_i' $ is connected for every $1 \le i \le k$.
\item[(F3)] Every $v \in V(G')\setminus T'$ has degree $0$ or $2$ in $G_i'$, for every $1 \le i \le k$.
\item[(F4)] Every $v \in T'$ satisfies $\sum_{i=1}^k {\mathit deg}_{G_i}(v) \le cap(v)+2$.
\end{itemize}
\end{lemma}

As in the Steiner tree case, we define $H$ to be the hypergraph on vertex set $T$ obtained from $G$ by replacing each vertex $v \in V(G) \setminus T$ with an edge of size 3 containing $N(v)$.  Similarly, we let $H'$ denote the hypergraph obtained from $G'$ by replacing each $v \in V(G') \setminus T'$ with an edge of size 3 containing $N(v)$.
For simplicity of notation, let $E=E(H)$ and $E'=E(H')$.
By this definition and equation~(\ref{gpdegconn}), every $v \in T'$ satisfies
\begin{equation}
\label{hpdegconn}
{\mathit deg}_{H'}(v) = g(k) - p(v) - q(v) - u(v).
\end{equation}
It follows from the connectivity properties of $G$ that the hypergraph $H$ will be $g(k)$-edge-connected, so $H'$ will not have too many small edge-cuts.
Define $p=\sum_{v\in T'}p(v)$.

\begin{lemma}
\label{edge_countconn}
If ${\mathcal P}'$ is a partition of $T'$ then
\[ 3 \lambda^{out}_{\mathcal P'}(E') \ge g(k) (|{\mathcal P'}| - 1) - g(k) + 2k+p \]
and if $|{\mathcal P}'| = 2$ this may be improved by $k + \tfrac{1}{2}p$ to
\[  3 \lambda^{out}_{\mathcal P'}(E') \ge 3k + \tfrac{3}{2}p. \]
\end{lemma}

\noindent{\it Proof:}
Extend ${\mathcal P}'$ to a partition ${\mathcal P}$ of $T$ by adding $\{w\}$ as a block.
Then ${\mathcal P}$ is a partition of $V(H)$ and since $H$ is $g(k)$-edge-connected, we have that every $B \in {\mathcal P}$ satisfies $d_H(B)=|\delta_H(B)| \ge g(k)$.
Let $W^{\circ}$ be the set of hyperedges in $H$ associated with an edge in $W^{\circ}_n \cup W^{\circ}_t$ (so $|W^{\circ}|=2k$ by (Q1)).
For $i=2,3$ let $W^{\circ}_i$ denote the set of hyperedges in
$W^{\circ}$ that have ends in exactly $i$ different blocks of ${\mathcal P}$.
Similarly, for $i=2,3$, we define $E_i$ to be the set of all edges not in $W^{\circ}$ which have ends in exactly $i$ different blocks of ${\mathcal P}$. In particular, $W^{\bullet}_t \subseteq E_2$ and $|W^{\bullet}_t|=p$, so $|E_2|\ge p$.
We now have
\begin{align*}
	3 \lambda^{out}_{\mathcal P'}(E')	
	&= 3 \lambda^{out}_{\mathcal P}(E)  + 3 |W^{\circ}_3 | - 3 g(k)\\
	&= 	\left(	 \sum_{B \in {\mathcal P}} d_{H}(B) + |E_2| + |W_2^{\circ}|\right)	+ 3 |W^{\circ}_3 | - 3g(k)\\
	&\ge	g(k) |{\mathcal P}| + p + |W_2^{\circ}| + |W^{\circ}_3 | - 3 g(k).
	\end{align*}	
Note in the second line that $\sum_{B\in \mathcal{P}}d_H(B)$ counts hyperedges with ends in 3 blocks of $\mathcal{P}$ three times, and those with ends in 2 blocks of $\mathcal{P}$ only twice.
Now, using $|W_2^{\circ}|+|W_3^{\circ}|=2k$ (by (Q1)) we obtain the desired inequality.

If $|{\mathcal P}'| = 2$ then we set $F$ to be the set of edges in $H$ which do not contain $w$
but do intersect both blocks of ${\mathcal P}'$.
Note that $W^{\bullet}_t \subseteq E_2\setminus F$ so $|E_2|\ge |F|+|W_t^{\bullet}|=|F|+p$.
We now repeat the same calculation.
\begin{align*}
	3 \lambda^{out}_{\mathcal P'}(E')	
	&= 3 \lambda^{out}_{\mathcal P}(E)  + 3 |W^{\circ}_3 | - 3 g(k)	\\
	&= 		\left( \sum_{B \in {\mathcal P}} d_{H}(B) + |E_2|  +  |W_2^{\circ}|\right)	+ 3 |W^{\circ}_3 | - 3g(k)	\\
	&\ge		3 g(k) + |F| + p + |W_2^{\circ}|  + 3|W^{\circ}_3 | - 3 g(k)	\\
	&\ge		\Big( |F| + |W^{\circ}_3| \Big) +  \Big( |W^{\circ}_2| + |W^{\circ}_3| \Big) +p	 \\
	&=		\lambda^{out}_{\mathcal P'} (E' ) + 2k  +p
\end{align*}
Rearranging the final inequality we get the desired result.
\hfill$\Box$

With this, we are finally ready to prove Lemma \ref{connector_lem}.

\bigskip

\noindent{\it Proof of Lemma \ref{connector_lem}:} We will establish this lemma by contradicting Lemma \ref{primeconn}, thus contradicting the assumption that $G$ is a counterexample.  We let $M$ be the hypergraphic matroid on $H'$, and let $r$ be the rank function of $M$.

 We claim that there cannot exist a vector $x \in {\mathbb R}_+^{E'}$ satisfying all of the following properties:
\begin{itemize}
\item $x(e) \le 1$ for every $e \in E'$
\item $x( \Lambda^{in}_{\mathcal P}( E') )  \le k r( {\mathcal P} )$ for every partition ${\mathcal P}$ of $T'$
\item $-x(E') \le - k r(M)$
\item $x( \delta_{H'}(v) ) \le {\mathit cap}(v)$ for every $v \in T'$.
\end{itemize}
To establish the claim, suppose, for a contradiction, that such an $x$ does exist.
Lemma \ref{partition_rank} tells us that $x$ is a fractional basis of $kM$ and $kM$ has rank $k\,r(M)$.
We apply Theorem \ref{kls} to the matroid $kM$, the vector $x$, and the family of subsets of $E'$ given by $\{ \delta_{H'} (v) \mid v \in T' \}$.
Since every edge is contained in at most three of these sets, the theorem gives us a basis $B$ of $kM$ for which
\begin{equation}\label{C3conn}
|B \cap \delta_{H'}(v)| \le cap(v)+2
\end{equation}
holds for every $v \in T'$.  Since $B$ is a basis of $kM$ we may choose a partition of $B$ into bases $\{B_1, \ldots B_k \}$ of $M$.
Fix $i \in \{1,\ldots,k\}$. By Theorem~\ref{thm22} we may replace every $e \in B_i$ with an edge $uv$, where $u,v \in e$, so that the resulting graph $F_i$ is a forest. Since $|B_i|=r(M)=|T'|-1$, $F_i$ is in fact a tree spanning $T'$. Now partition $B_i$ into $B_i^2$ and $B_i^3$, where $B_i^j$ is the set of hyperedges in $B_i$ of size $j$.   Then $B_i^2 \subseteq E(F_i)$ and every $e \in B_i^3$ corresponds to an edge in $F_i$. Recall that $e=\{u_1,u_2,u_3\} \in B_i^3$ corresponds to a vertex $v \in V(G')\setminus T'$ with $N_{G'}(v)=\{u_1,u_2,u_3\}$. We construct a $T'$-connector $G_i$ from $B_i$ as follows. Start with $G_i'=(T',B_i^2)$. For every $e=\{u_1,u_2,u_3\} \in B_i^3$, add to $G_i'$ the vertex $v$ and the edges $u_1v$ and $u_2v$, where $u_1u_2$ is the edge in $F_i$ corresponding to $e$, and $v$ is the vertex in $V(G')\setminus T'$ with $N_{G'}(v)=\{u_1,u_2,u_3\}$.
Since $F_i$ is a tree spanning $T'$, $G_i'$ is connected and spans $T'$. Moreover every vertex in $V(G')\setminus T'$ has degree $0$ (if $N_{G'}(v) \notin B_i$) or $2$ (if $N_{G'}(v) \in B_i$) in $G_i'$.
Therefore each $G_i'$ satisfies (F1), (F2), (F3) and (F4). Thus, these graphs contradict  Lemma \ref{primeconn}, and our claim has been established.

Since we know that $x$ cannot exist, Farkas' Lemma tells us that there must exist families of nonnegative parameters $\{y_e\}_{e \in E'}$, $y_G$, $\{ y_v \}_{v \in T'}$, and $\{ y_{\mathcal P} \}$ for every partition ${\mathcal P}$ of $T'$, that satisfy the following properties
\begin{equation}
\label{eqa1conn}
y_G\, k (|T'|-1) > \sum_{e \in E'} y_e + \sum_{ {\mathcal P} } y_{\mathcal P} \, k \, r({\mathcal P}) + \sum_{v \in T'} y_v \, {\mathit cap}(v)
\end{equation}
\begin{equation}
\label{eqb1conn}
y_e + \sum_{ {\mathcal P} : e \in \Lambda^{in}_{\mathcal P}(E') } y_{\mathcal P} + \sum_{v \sim e} y_v \ge y_G \mbox{\quad\quad for every $e \in E'$}
\end{equation}

As in the Steiner tree case, we may assume that $y_G=1$, that for the trivial partition $\{T'\}$, $y_{ \{T'\} } = 0$, and that the only partitions ${\mathcal P}$ with $y_{ \mathcal P} > 0$ are nested.
Since the trivial partition $\{T'\}$ has zero weight, we may assume that there are partitions ${\mathcal P}_0 \le {\mathcal P}_1 \le \ldots \le {\mathcal P}_{\ell}$ so that every partition ${\mathcal Q} \not\in \{ {\mathcal P}_0, \ldots, {\mathcal P}_{\ell} \}$ has $y_{\mathcal Q} = 0$ and so that
${\mathcal P}_0$ is the discrete partition and ${\mathcal P}_{\ell}$ is a partition of size two.
For convenience of notation, we now set $y_i = y_{ {\mathcal P}_i }$, $\Lambda^{out}_i(F) = \Lambda^{out}_{ {\mathcal P}_i }(F)$, $\Lambda^{in}_i(F) = \Lambda^{in}_{ {\mathcal P}_i }(F)$ and assign $\lambda_i^{out}$ and $\lambda^{in}_i$ similarly for every $0 \le i \le \ell$ and every subset of edges $F$.  With this, we now restate equation (\ref{eqb1conn}) as follows:
\begin{equation}
\label{eqb2conn}
y_e + \sum_{ i \in \{0,\ldots, \ell\} : e \in \Lambda_i^{in}(E') } y_i + \sum_{v \sim e} y_v \ge 1 \mbox{ \hspace*{.2in} for every $e \in E'$}.
\end{equation}

As in Section~\ref{sec:steiner}, we may assume that $\sum_{i=0}^{\ell} y_i = 1$.
We use this normalization next to rewrite the following term:
\begin{align*}
&k( |T'| - 1 ) - \sum_{i=0}^{\ell} y_i k r( {\mathcal P}_i) = \sum_{i=0}^{\ell} y_i k \big( |T'| - 1 - r ( {\mathcal P_i} ) \big) = \sum_{i=0}^{\ell} y_i k ( |{\mathcal P}_i| - 1 ).
\end{align*}
Taking advantage of this identity we may now rewrite equation (\ref{eqa1conn}) as follows
\begin{equation}
\label{eqa2conn}
\sum_{i=0}^{\ell}  y_i \, k( |{\mathcal P}_i| - 1)  > \sum_{e \in E'} y_e + \sum_{v \in T'} y_v \, {\mathit cap}(v).
\end{equation}

In the remainder of our proof we will use (\ref{eqa2conn}) and (\ref{eqb2conn}), along with some previous results, to achieve our desired contradiction. In the following calculation, the first inequality is obtained by summing (\ref{eqb2conn}) over all edges, the second inequality uses Lemma \ref{edge_countconn} (including the special case that $|{\mathcal P}_{\ell}| = 2$), and the final inequality uses equation (\ref{eqa2conn}):
\begin{align*}
 3 \sum_{e \in E'} y_e &+ 3\sum_{v \in T'} y_v {\mathit deg_{H'}}(v)
	\ge 3|E'| -  3\sum_{j=0}^{\ell} y_j  \, \lambda^{in}_j (E') 	\\
	&= 3\sum_{j=0}^{\ell}  y_{j} \, \lambda^{out}_j(E') 	\\
	&\ge \sum_{j=0}^{\ell}   y_j \Big( g(k)( |{\mathcal P}_j| - 1) 		
			- g(k) + 2k +p \Big)   + y_{\ell} (k+\tfrac{1}{2}p) \\
	&\ge \tfrac{g(k)}{k}  \sum_{j=0}^{\ell}   y_j \, k( |{\mathcal P}_j| - 1)		
			-  g(k) + 2k + p+ y_{\ell} k\\
	&> \tfrac{g(k)}{k}\Big( \sum_{e \in E'} y_e + \sum_{v \in T'} y_v  {\mathit cap}(v)\Big)
		- g(k) + 2k +p+ y_{\ell} k.
\end{align*}

Rearranging we obtain
\begin{equation}\label{eqfk2conn0} g(k) - 2k > p+ \tfrac{1}{k} \sum_{v \in T'} y_v \Big( g(k) {\mathit cap}(v) - 3k {\mathit deg}_{H'}(v) \Big) + \tfrac{g(k) - 3k}{k}  \sum_{e \in E'} y_e
	+ y_{\ell} k .\end{equation}
	
Equation~(\ref{eqb2conn}) implies that $y_v$ never needs to be larger than 1, therefore we may assume that $y_v\leq 1$ for every $v\in T'$ and so
$p\geq \sum_{v\in T'}p(v)y_v$.
With this fact and equation~(\ref{eqfk2conn0}) we obtain
\begin{equation}\label{eqfk2conn} g(k) - 2k > \tfrac{1}{k} \sum_{v \in T'} y_v \Big( g(k) {\mathit cap}(v) - 3k {\mathit deg}_{H'}(v) +kp(v)\Big) + \tfrac{g(k) - 3k}{k}  \sum_{e \in E'} y_e
	+ y_{\ell} k .\end{equation}
	

In the next calculation we consider an arbitrary vertex $v \in T'$.  We call upon the definition of the capacity function, ${\mathit cap}(v) = g(k) - 2k - 2 - p(v) - q(v)$ (equation (\ref{capdefconn})) and the degree formula ${\mathit deg}_{H'}(v) = g(k) - p(v) -q(v) - u(v)$ (equation (\ref{hpdegconn})); the inequalities $q(v) \le u(v) + 1$ (Observation \ref{obp2uconn}), $p(v) \le k-1$ (Lemma \ref{buzzer_beaterconn}) and $p(v)+q(v)\leq 3k+3$ (Lemma~\ref{p1plusp2}); as well as the value of $g(k)=6k+6$.
\begin{align*}
&g(k) {\mathit cap}(v) - 3k {\mathit deg}_{H'}(v)+kp(v)\\
	&= g(k) (4k+4-p(v) - q(v))  - 3k (g(k) - p(v) - q(v) - u(v))+kp(v)	\\
	&= (6k+6)(4k+4-p(v)-q(v))-3k(6k+6)+4kp(v)+3kq(v)+3ku(v)\\
	&= (6k+6) (k+4) - (2k+6)p(v)-6q(v)-3k(q(v)-u(v))	\\
	&\geq  (6k+6) (k+4) - (2k+6)p(v)-6q(v)-3k\\
	&= (6k+6) (k+4)-3k - 2kp(v)-6(p(v)+q(v))\\
	&\ge (6k+6) (k+4)-3k - 2k(k-1)-6(3k+3)\\
	&=4k^2+11k+6\\
	&\ge k(g(k) - 2k).
\end{align*}
Combining this with (\ref{eqfk2conn}) gives us the following:
\begin{equation}\label{fk2kconn} (g(k) - 2k) \Big( 1 - \sum_{v \in T'} y_v \Big) > \tfrac{g(k) - 3k}{k} \sum_{e \in E'} y_e  + y_{\ell}k. \end{equation}

Setting $\alpha = 1 - \sum_{v \in T'} y_v$ it follows from (\ref{fk2kconn}) that $\alpha > 0$.  Furthermore, every edge $e \in \Lambda^{out}_{\ell}(E')$ must satisfy $y_e \ge \alpha$ by (\ref{eqb2conn}), since this edge is outer with respect to every partition ${\mathcal P}_i$ and $\sum_{v \sim e} y_v \le \sum_{v \in T'} y_v = 1- \alpha$.  We are now able to derive a contradiction based on the value of $y_{\ell}$.

First suppose that $y_{\ell} \ge \alpha$.  In this case (\ref{fk2kconn}) gives us the following contradiction (here we have applied Lemma \ref{edge_countconn} for the last inequality).
\begin{align*}
(g(k) - 2k) \alpha
	&> \tfrac{g(k) - 3k}{k} \sum_{e \in E'} y_e + y_{\ell} k	\\
	&\ge \tfrac{g(k)-3k}{k} \alpha \,\lambda^{out}_{\ell}(E' )  + \alpha k	\\
	&\ge \tfrac{g(k)-3k}{k} \alpha k + \alpha k	\\
	&= (g(k) - 2k) \alpha.
\end{align*}

Next assume that $y_{\ell} \leq \alpha$. Then  $\beta = \alpha - y_{\ell} \ge 0$.  Now, every edge in $\Lambda^{out}_{\ell-1}(E') \setminus \Lambda^{out}_{\ell}(E' )$ must satisfy
$y_e \ge \beta$  by (\ref{eqb2conn}), since this edge receives a total contribution of at most $1- \alpha$ from the vertices and $y_{\ell}$ from the partitions.  Hence (\ref{fk2kconn}) gives us
\[
(g(k) - 2k) \alpha >  \tfrac{g(k) - 3k}{k} (\alpha - \beta) \lambda^{out}_{\ell}(E') + \tfrac{g(k)-3k}{k} \beta \, \lambda^{out}_{\ell-1}(E') + (\alpha - \beta) k.
\]
Splitting the left hand side into $(g(k) - 2k) (\alpha - \beta) + (g(k) - 2k) \beta$ and then rearranging gives:
\begin{equation}\label{alphabeta} 0 > (\alpha - \beta) \Big( \tfrac{g(k) - 3k}{k} \lambda^{out}_{\ell}(E' ) + k  - (g(k) - 2k) \Big)
	+ \beta \Big( \tfrac{g(k)-3k}{k} \lambda^{out}_{\ell-1}(E')  - (g(k) - 2k) \Big). \end{equation}
As in the previous case the term above with coefficient $\alpha - \beta$ must be $\ge 0$ by Lemma \ref{edge_countconn}.  Using the same lemma, we find that
$$\lambda^{out}_{\ell-1}(E') \ge \tfrac{1}{3} (g(k) + 2k) \ge 2k$$
and since the coefficient of $\beta$ in (\ref{alphabeta}) must be negative we have
\[ 0 > \big( \tfrac{g(k)-3k}{k}(2k) - ( g(k) - 2k ) \big) = g(k) - 4k \ge 0,  \]
and this final contradiction completes the proof.
\hfill$\Box$

\section{Conclusion}

Our methods leave some room for improvement in the values of $f(k)$ and $g(k)$, but not by more than about $k$. To see this, consider a hypergraph $R$ where all edges are of size three and where $E(R)$ is an independent set in the associated hypergraphic matroid. The average degree of a vertex in the hypergraph is $\tfrac{3|E(R)|}{|V(R)|}=\tfrac{3(|V(R)|-1)}{|V(R)|}=3-\tfrac{3}{|V(R)|} \sim 3$.
Therefore, in our proof, when we look for a fractional basis $x$ of a union of $k$ hypergraphic matroids  and examine a vertex $v$, we cannot hope that
$$x(\delta(v))+ p(v)+ q(v)$$
is really any smaller than $3k$. However, in our LPs, we require
$$x(\delta(v)) \leq f(k)-k-2-p(v)-q(v)$$
and
$$x(\delta(v)) \leq g(k)-2k-2-p(v)-q(v),$$
respectively. Hence we will have to have
$$3k-p(v)-q(v) \leq f(k)-k-2-p(v)-q(v),$$
and
$$3k-p(v)-q(v) \leq g(k)-2k-2-p(v)-q(v),$$
which implies that $f(k) \geq 4k+2$ and  $g(k)\geq 5k+2.$ This would seem to be the limit of our current approach. Moreover, the fight for this extra $k$ appears to require a great deal more case analysis than we are willing to do at this point.

\section*{Acknowledgement}
The authors would like to thank the many adorable kittens who gave their lives to make this research possible.  We will always love you Biscuit, Paco, Sunshine, Pookie, French Fry, Goober, Chicken Noodle, Buttons, Pancake, Snowball, Peaches, and friends.


\begin{thebibliography}{00}

\bibitem{Edmonds} J. Edmonds. Submodular functions, matroids and certain polyhedra, Combinatorial structures and their applications (Proc. Calgary Internat. Conf. 1969), 69--87. Gordon and Breach, New York.

\bibitem{ER} J. Edmonds and G.-C. Rota. Submodular set functions (Abstract). Waterloo Combinatorics conference (1966).


\bibitem{FKK} A. Frank, T. Kir\'{a}ly, and M. Kriesell. On decomposing a hypergraph into k connected
sub-hypergraphs, Discrete Appl. Math. 131 (2003), 373--383.


\bibitem{La} L. C. Lau. An approximate max-Steiner-tree-packing min-Steiner-cut theorem, Combinatorica 27 (2007), 71--90.

\bibitem{Lorea} M. Lorea. Hypergraphes et matroides, Cahiers Centre Etud. Rech. Oper. 17 (1975), 289--291.

\bibitem{KLS} T. Kir\'{a}ly, L. C. Lau, and M. Singh. Degree bounded matroids and submodular flows,
Proceedings of 13th International Conference IPCO 2008, LNCS 5035 (2008), 259--272.

\bibitem{Kr1} M. Kriesell. Edge-disjoint trees containing some given vertices in a graph, J. Combin.
Theory Ser. B 88 (2003), 53--65.

\bibitem{Kr2 }M. Kriesell. Edge disjoint Steiner trees in graphs without large bridges, J. Graph Theory
62 (2009), 188--198.

\bibitem{Ma} W. Mader. A reduction method for edge-connectivity in graphs, Ann. Discrete Math. 3
(1978), 145--164.

\bibitem{NW} C. St. J. A. Nash-Williams, Edge-disjoint spanning trees of finite graphs, J. London Math. Soc. 36 (1961), 445--450.

\bibitem{Rado} R. Rado. A theorem on independence relations, Quart. J. Math. Oxford 13 (1942), 83--89.

\bibitem{Tu} W. T. Tutte. On the problem of decomposing a graph into n connected factors, J. London Math. Soc. 36 (1961), 221--230.

\bibitem{WW} H. Wu and D. B. West. Packing of Steiner trees and S-connectors in graphs. J. Combin. Theory Ser. B, 102 (2012), 186--205.

\end{thebibliography}
\end{document}